\documentclass[10pt]{amsart}
\usepackage{amsfonts,amsthm,amsmath,color,graphics,hyperref}
\usepackage[all]{xy}
\usepackage{indentfirst}
\newtheorem{Theorem}{Theorem}[section]
\newtheorem{Corollary}[Theorem]{Corollary}
\newtheorem{Definition}[Theorem]{Definition}
\newtheorem{Lemma}[Theorem]{Lemma}
\newtheorem{Proposition}[Theorem]{Proposition}
\newtheorem{Remark}[Theorem]{Remark}

\newtheorem{Example}[Theorem]{Example}

\begin{document}
\title{Weierstrass filtration on Teichm\"{u}ller curves and Lyapunov exponents}
\author{Fei Yu}
\author{Kang Zuo}
\address{School of Mathematical Sciences, Xiamen University, Xiamen, 361005, People¡¯s Republic of China}
\address{Fachbereich 08-Physik Mathematik und Informatik, Universit{\"a}t Mainz, 55099 Mainz, Germany}
\email{vvyufei@gmail.com,yu@xmu.edu.cn}
\address{Fachbereich 08-Physik Mathematik und Informatik, Universit{\"a}t Mainz, 55099 Mainz, Germany}
\email{zuok@uni-mainz.de}

 \maketitle
\begin{abstract}
We define the Weierstrass filtration for Teichm\"{u}ller curves and
construct the Harder-Narasimhan filtration of the Hodge bundle of a
Teichm\"{u}ller curve in hyperelliptic loci and low-genus
nonvarying strata. As a result we obtain the sum of Lyapunov
exponents of Teichm\"{u}ller curves in these strata.
\end{abstract}
\tableofcontents
 \footnote{


This work is supported by the SFB/TR 45 ¡®Periods,
Moduli Spaces and Arithmetic of Algebraic Varieties¡¯ of the DFG.

The first author is also supported by the Fundamental Research Funds
for the Central Universities  (No. 2013121001).}
\section{Introduction}
The computation of Lyapunov exponents is an important subject in the
theory of Teichm\"{u}ller curves. The sum of Lyapunov exponents for
hyperelliptic connected components is known by \cite{EKZ}. In
general it is determined by the Siegel-Veech constants which measure
the boundary behavior of the underlying $SL_2(\mathbb{R})$-orbit
closure, cf.\ \cite{EMZ03}\ \cite{CM11}\ \cite{KZ03}. In low-genus
cases \cite{Ba}\ \cite{LM11}\ \cite{Mo}, the Lyapunov spectrum has
been worked out in some special cases, based on a series of
work by McMullen \cite{Mc03}\ \cite{Mc06a}\ \cite{Mc07}. Some
concrete examples like square-tiled surfaces and triangle groups
were computed in \cite{BM10}\ \cite{EKZ11}\ \cite{FMZ11a}.

The sums of nonnegative Lyapunov exponents of Teichm\"{u}ller curves
in a number of strata components of low-genus flat surfaces are
constants. The same result holds for the hyperelliptic loci of any
stratum. The proofs of these non-varying results are completely
different \cite{Mo12}. For the first, one uses algebraic geometry,
in particular slope calculations \cite{CM11}\ \cite{Chen}, and the
other relies on the correspondence to Siegel-Veech constants
\cite{EKZ}.

In this paper, we will construct Weierstrass filtrations of the
Hodge bundle based on the dimension of sublinear systems of zeros of
holomorphic differentials, define Weierstrass exponents according to
the Harder-Narasimhan filtration, and compute them in the
hyperelliptic loci and low-genus nonvarying strata. This will allow
us to give a unified method to compute the sum of Lyapunov
exponents.

  Let $g\geq1$ be an integer, and let $(m_1,...,m_k)$ be a partition of $2g-2$. Denote by
$\Omega\mathcal{M}_g(m_1,...,m_k)$ the stratum parameterizing genus
$g$ Riemann surfaces with Abelian differentials that have $k$
distinct zeros of order $m_1,...,m_k$, respectively.

Let $C$ be a Teichm\"{u}ller curve that lies in
$\Omega\mathcal{M}_g(m_1,...,m_k)$. Denote by $f:S\rightarrow C$ the
universal family over a Teichm\"{u}ller curve with  distinct
sections $D_1,...,D_k$. The relative canonical bundle can be
computed through the formula \eqref{canonical}:
$$\omega_{S/C}\simeq f^*\mathcal{L}\otimes \mathcal{O}(\Sigma m_i D_i).$$

In  hyperelliptic loci and low-genus nonvarying strata, we will
construct the Harder-Narasimhan filtration of $f_*(\omega_{S/C})$
and show that the factors of the Jordan-H\"{o}lder filtration of
each semistable graded quotient are line bundles (see the filtration
\eqref{exponent} below). Write

$$0\subset V_1\subset V_2\subset...\subset V_g=f_*(\omega_{S/C})$$
for the filtration, then the $i$-th Weierstrass exponent $w_i$ is
defined as
$$\frac{\mathrm{deg}(V_i/V_{i-1})}{\mathrm{deg}(\mathcal{L})}$$
$(i=1,...,g$).
\begin{Theorem}[Theorem \ref{hyperelliptic}] Let $C$ be a Teichm\"{u}ller curve in the hyperelliptic locus of some
stratum $\overline{\Omega\mathcal{M}}_g(m_1,...,m_k)$, and denote by
$(d_1,...,d_n)$ the orders of singularities of underlying quadratic
differentials. Then the Weierstrass exponent $w_i$ for $C$ is the
$i$-th largest number in the following set
$$\{1\}\cup\Big\{1-\frac{2k}{d_j+2}\Big\}_{\forall d_j, 0<2k\leq d_j+1}.$$
\end{Theorem}
This result can be used to recover the sum of Lyapunov exponents in
hyperelliptic loci originally found in \cite{EKZ}:
$$\sum \lambda_i=\sum w_i=\frac{1}{4}\sum_{d_j\mathrm{\ odd} } \frac{1}{d_j+2}.$$ It was
conjectured by Kontsevich and Zorich in \cite{KZ97} (for
Teichm\"{u}ller geodesic flows), and was shown by M.Bainbridge
in the case $g=2$ \cite{Ba}.

Zorich communicated to D.Chen and M.M\"{o}ller that, based on a
limited number of computer experiments about a decade ago,
Kontsevich and Zorich observed that the sum of Lyapunov exponents is
nonvarying among all the Teichm\"{u}ller curves in a stratum
roughly if the genus plus the number of zeros is less than seven,
while the sum varies if this sum is greater than seven. The
following two results are entirely based on the paper \cite{CM11}.
They have verified those nonvarying strata observed by Kontsevich
and Zorich for $g\leq 5$, except for the strata
$\overline{\Omega\mathcal{M}}^{even}_4(4,2),\overline{\Omega\mathcal{M}}^{odd}_4(4,2)$
and $\overline{\Omega\mathcal{M}}^{odd}_5(6,2)$. The nonvarying
result is obtained by showing that Teichm\"{u}ller curves in a
stratum are disjoint with a geometrically defined divisor on the
moduli spaces of curves. Our main result provides another proof of
the nonvarying property which also solves the remaining case with
the help of D.Chen.

\begin{Theorem} [Theorem \ref{2g-2}] Let $k_i$ be the $i$-th
largest number in the Weierstrass gap $G_p$. For a Teichm\"{u}ller curve in
$\overline{\Omega\mathcal{M}}^{hyp}_g(2g-2)$, the Weierstrass
exponent $w_i$ is $\frac{k_i}{2g-1}$ and the Weierstrass gap $G_p$
equals $\{1, 3, 5,...,2g-5,2g-3,2g-1\}$.

If moreover $g\leq 5$, then
\begin{enumerate}
\item[1.] For a Teichm\"{u}ller curve in
$\overline{\Omega\mathcal{M}}^{odd}_g(2g-2)$, the Weierstrass
exponent $w_i$ is $\frac{k_i}{2g-1}$ and the Weierstrass gap $G_p$
equals $\{1, 2, 3,...,g-2,g-1,2g-1\}$, and $f_*{\omega_{S/C}}$
splits into direct sum of line bundles.
\item[2.] For a Teichm\"{u}ller curve in
$\overline{\Omega\mathcal{M}}^{even}_g(2g-2)$, the Weierstrass
exponent $w_i$ is $\frac{k_i}{2g-1}$ and the Weierstrass gap $G_p$
equals $\{1, 2, 3,...,g-2,g,2g-1\}$.
\end{enumerate}
\end{Theorem}

\begin{Theorem}[Theorem \ref{lowgenus}] For a Teichm\"{u}ller curve in the strata\\
$$\overline{\Omega\mathcal{M}}_3(3,1),\overline{\Omega\mathcal{M}}^{odd}_3(2,2),\overline{\Omega\mathcal{M}}_3(2,1,1)$$
$$\overline{\Omega\mathcal{M}}_4(5,1),\overline{\Omega\mathcal{M}}^{odd}_4(4,2),\overline{\Omega\mathcal{M}}^{non-hyp}_4(3,3),\overline{\Omega\mathcal{M}}^{odd}_4(2,2,2),\overline{\Omega\mathcal{M}}_4(3,2,1)$$
$$\overline{\Omega\mathcal{M}}_5(5,3),\overline{\Omega\mathcal{M}}^{odd}_5(6,2),$$
the Weierstrass exponents can explicitly be calculated as in Tables
1, 2 and 3. Moreover $f_*{\omega_{S/C}}$ splits into a direct sum of
line bundles.

For a Teichm\"{u}ller curve in the stratum
$\overline{\Omega\mathcal{M}}^{even}_4(4,2)$, the Weierstrass
exponents can explicitly be calculated as in Table 2.
\end{Theorem}
Related work about quadratic differentials was done in
\cite{CM12}.

Our basic idea is to construct filtrations of
$f_*\mathcal{O}(\omega_{S/C})$
 $$0\subset
\mathcal{L}\subset...\subset f_*\mathcal{O}\Big(\omega_{S/C}-\sum
d_iD_i\Big) \subset...\subset f_*\mathcal{O}(\omega_{S/C})$$ 

and
then compute each graded quotient. But generally, it is difficult to
compute the quotient. The quotient is a locally free subsheaf of the
direct image of a bundle $\mathcal{O}_{aD_i}(dD_i)$, hence we use
the Harder-Narasimhan filtration to study the bundle. The difficulty
will disappear if we assume that the Weierstrass semigroup of fibers
is non-varying.

The paper is organized as follows. In section 2, we introduce some
background material that has appeared in \cite{CM11}. In section 3,
we give a basic example to construct the Weierstrass filtration and
the Weierstrass semigroup filtration. In some special cases, we also
show that it is the Harder-Narasimhan filtration, and compute
Weierstrass exponents of such Teichm\"{u}ller curves. In section 4
we define Weierstrass filtrations and apply it to compute the sum of
Lyapunov exponents. In section 5, we define We
ierstrass exponents
and compute them for the non-varying strata.

\section{Background}
\subsection{Moduli spaces}

Denote by $\Omega\mathcal{M}_g$  the moduli space of pairs
$(X,\omega)$ where $X$ is a curve of genus $g$ and $\omega$ is a
holomorphic one-form on $X$. It is fibred over the moduli space
$\mathcal{M}_g$ of curves. Let $(m_1,...,m_k)$ be a partition of
$2g-2$, and let $\Omega\mathcal{M}_g(m_1,...,m_k)$ denote the
stratum parameterizing one-forms that have $k$ distinct zeros of
order $m_1,...,m_k$, respectively.  Denote by
$\Omega\mathcal{M}^{hyp}_g(m_1,...,m_k)$( resp. odd, resp. even) the
hyperelliptic (resp.  odd theta characteristic, resp. even theta
characteristic) connected component (\cite{KZ03}).

Let $\overline{\mathcal{M}}_g$ denote the Deligne-Mumford
compactification of $\mathcal{M}_g$. Then $\Omega\mathcal{M}_g$
extends as a vector bundle over $\overline{\mathcal{M}}_g$,
parameterizing sections of the dualizing sheaf or equivalently
stable one-forms. We denote by $\overline{\Omega\mathcal{M}}_g$ the
total space of this extension. Points in
$\overline{\Omega\mathcal{M}}_g$, called flat surfaces, are also
written as $(X,\omega)$ with $\omega$ a stable one-form on $X$.

 Let $C$ be a genus-$g$ curve and $\mathcal{L}$ a line bundle of
degree $d$ on $C$. Denote by $|\mathcal{L}|$ the projective space of
one-dimensional subspaces of $H^0(C,\mathcal{L})$. For a
(projective) $r$-dimension linear subspace $V$ of $|\mathcal{L}|$,
we call $(\mathcal{L},V)$ a linear series of type $g^r_d$.

Let $\underline{\omega}=(\omega_1,...,\omega_n)$ be a tuple of
integers. The generalized Brill-Noether locus
$BN^r_{d,\underline{w}}$ is the locus in
$\overline{\mathcal{M}}_{g,n}$ of pointed curves $(C,p_1,...,p_n)$
with a line bundle $\mathcal{L}$ of degree $d$ such that
$\mathcal{L}$ admits a linear system $g^r_d$  and
$h^0(\mathcal{L}(-\sum w_ip_i))\geq 1$.

We need the following generalization of Clifford's Theorem for
stable curves:
\begin{Theorem}[{\cite[p. 8, Theorem 2.5]{CM11}}] \label{clifford} Let $C$ be a stable curve and $D$ an effective divisor
with $\mathrm{deg}(D)\leq 2g-1$. Then
$$h^0(\mathcal{O}_C(D))-1\leq \mathrm{deg}(D)/2,$$
if one of the following conditions holds:
\begin{enumerate}
\item[1.] $C$ is smooth;
\item[2.] $C$ has at most two components;
\item[3.] $C$ does not have separating nodes and $\mathrm{deg}(D)\leq 4$.
\end{enumerate}
\end{Theorem}
\subsection{Teichm\"{u}ller curves}
A Teichm\"{u}ller curve $C$ is an algebraic curve in $\mathcal{M}_g$
that is totally geodesic with respect to the Teichm\"{u}ller metric.
After suitable base change, we can get a universal family
$f:S\rightarrow C$, which is a relatively minimal semistable model
with disjoint sections $D_1,...,D_k$; here $D_i|_X$ is a zero of
$\omega$ when restricting to each fiber $X$ (\cite[p. 11]{CM11}).

Let $\mathcal{L}\subset f_*{\omega_{S/C}}$ be the line bundle whose
fiber over the point in $C$ corresponding to $X$ is
$\mathbb{C}\omega$, the generating differential of $C$; it is also
known as  the ''maximal Higgs'' line bundle, in the sense of
\cite{VZ} and \cite{Mo11}. Let $\Delta\subset \overline{B}$ be the
set of points with singular fibers. Then the property of being
''maximal Higgs'' says by definition that $\mathcal{L}\cong
\mathcal{L}^{-1}\otimes\omega_C(log{\Delta})$ and
$$\mathrm{deg}(\mathcal{L})=(2g(C)-2+|\Delta|)/2,$$
together with an identification (relative canonical bundle formula
\cite[p. 19]{CM11}):
\begin{equation}\label{canonical}
 \omega_{S/C}\simeq f^*\mathcal{L}\otimes \mathcal{O}(\Sigma m_i
 D_i).
\end{equation}
By the adjunction formula we get
$$D^2_i=-\omega_{S/C}D_i=-m_iD^2_i-\mathrm{deg}{\mathcal{L}}$$
and thus
\begin{equation}\label{intersection}
 D^2_i=-\frac{1}{m_i+1}\mathrm{deg}\mathcal{L}.
\end{equation}

The variation of Hodge structures (VHS for short) over a
Teichm\"{u}ller curve decomposes into sub-VHS
\begin{equation}\label{VHS}
R^1f_*\mathbb{C}=(\bigoplus^r_{i=1} \mathbb{L}_i)\oplus \mathbb{M}.
\end{equation}
Here $\mathbb{L}_i$ are rank-2 subsystems, maximal Higgs
$\mathbb{L}^{1,0}_1\simeq \mathcal{L}$ for $i=1$, nonunitary but
not maximal Higgs  for $i\neq 1$ \cite[Theorem 2.2]{Mo11}.

Here we collect some properties of Teichm\"{u}ller curves along the
boundary of the moduli space which will be needed  in the subsequent sections.
\begin{Theorem}[{\cite[section 3.3]{CM11}}]\label{singular}
\begin{enumerate}
\item[1.] The section $\omega$ of the canonical bundle of each smooth
fiber over a Teichm\"{u}ller curve extends to a section
$\omega_{\infty}$ for each singular fiber $X_{\infty}$ over the
closure of a Teichm\"{u}ller curve. The signature of zeros of
$\omega_{\infty}$ is the same as that of \ $\omega$. Moreover,
$X_{\infty}$ does not have separating nodes.
\item[2.] For Teichm\"{u}ller curves generated by a flat surface in
$\Omega\mathcal{M}_g(2g-2)$ the degenerating fibers are irreducible.
\item[3.] Let $C$ be a Teichm\"{u}ller curve generated by an Abelian
differential $(X,\omega)$ in $\Omega\mathcal{M}_g(\mu)$. Suppose
that an irreducible degenerating fiber $X_\infty$ over a cusp of $C$
is hyperelliptic. Then $X$ is hyperelliptic, hence the whole
Teichm\"{u}ller curve lies in the locus of hyperelliptic flat
surfaces.
\item[4.] Let $C$ be a Teichm\"{u}ller curve generated by a flat surface
in $\Omega\mathcal{M}_5(8)$ even. Then $C$ does not intersect the
Brill-Noether divisor $BN^1_3$ on $\mathcal{M}_5$.
\item[5.] Moreover, if \ $\mu$ lies in $ \{(4), (3, 1), (6), (5, 1), (3,
3), (3, 2, 1), (8), (5, 3)\}$ and $(X,\omega)$ is not hyperelliptic,
then the nondegenerating fibers of the Teichm\"{u}ller curve are
hyperelliptic.
\end{enumerate}
\end{Theorem}

\subsection{Lyapunov exponents}
Fix an $SL_2(\mathbb{R})$-invariant, ergodic measure $\mu$ on
$\Omega\mathcal{M}_g$. Let $V$ be the restriction of the real Hodge
bundle (i.e. the bundle with fibers $H^1(X,\mathbb{R})$) to the
support $M$ of $\mu$. Let $S_t$ be the lift of the geodesic flow to
$V$ via the Gauss-Manin connection. Then Oseledec's multiplicative
ergodic Theorem guarantees the existence of a filtration
$$0\subset V_{\lambda_g}\subset ...\subset V_{\lambda_1}=V$$
by measurable vector subbundles with the property that, for almost
all $m\in M$ and all $v\in V_m\backslash\{0\}$ one has
$$||S_t(v)||=\mathrm{exp}(\lambda_it+o(t)),$$
where $i$ is the maximal index such that $v$ is in the fiber of
$V_i$ over $m$ i.e.$v\in(V_i)_m$. The numbers $\lambda_i$ for
$i=1,...,k\leq rank(V)$ are called the \emph{Lyapunov exponents} of
$S_t$. Since $V$ is symplectic, the spectrum is symmetric in the
sense that $\lambda_{g+k}=-\lambda_{g-k+1}$. Moreover, from
elementary geometric arguments it follows that one always has
$\lambda_1=1$. Thus, the Lyapunov spectrum is defined by the
remaining non-negative Lyapunov exponents
$$\lambda_2\geq...\geq\lambda_g.$$
The bridge between the 'dynamical' definition of Lyapunov exponents
and the  'algebraic' method applied in the sequel is given by the
following result.

\begin{Theorem}[{\cite{KZ97} \cite{BM10}}]\label{sumly} If the VHS over the Teichm\"uller curve $C$ contains a sub-VHS $\mathbb{W}$
of rank $2k$, then the sum of the $k$ corresponding non-negative
Lyapunov exponents equals
$$\overset{k}{\underset{i=1}{\sum}}\lambda^{\mathbb{W}}_i=\frac{2\mathrm{deg} \mathbb{W}^{(1,0)}}{2g(C)-2+|\Delta|},$$
where $\mathbb{W}^{(1,0)}$ is the $(1,0)$-part of the Hodge
filtration of the vector bundle associated with $\mathbb{W}$. In
particular, we have
$$\overset{g}{\underset{i=1}{\sum}}\lambda_i=\frac{2\mathrm{deg}f_*\omega_{S/C}}{2g(C)-2+|\Delta|}.$$
\end{Theorem}
Let $L(C)=\overset{g}{\underset{i=1}{\sum}}\lambda_i$ be the sum of
Lyapunov exponents, and define
$$\kappa_{\mu}=\frac{1}{12}\overset{k}{\underset{i=1}{\sum}}\frac{m_i(m_i+2)}{m_i+1}.$$
Eskin, Kontsevich and Zorich have a formula to compute $L(C)$ (for
the Teichm\"{u}ller geodesic flow):
\begin{Theorem}[{\cite[Theorem 1]{EKZ}}] For the VHS over the Teichm\"uller curve $C$, we have
$$L(C)=\kappa_{\mu}+\frac{\pi^2}{3}c_{area}(C),$$
where $c_{area}(C)$ is the area Siegel-Veech constant of $C$.
\end{Theorem}

\subsection{Vector bundles on curves}The readers can refer to \cite[section 1.3, section 1.5]{HL}
for details about the Harder-Narasimhan filtration and the
Jordan-H\"{o}lder filtration of sheaves on a variety. Let $C$ be a
smooth curve and $V$ a vector bundle over $C$ of slope
$\mu(V):=\frac{\mathrm{deg}(V)}{\mathrm{rk}(V)}$. We call $V$
semistable (resp. stable) if $\mu(W)\leq\mu(V)$ (resp.
$\mu(W)<\mu(V)$) for any subbundle $W\subset V$.

A Harder-Narasimhan filtration for $V$ is an increasing filtration:
$$0=HN_0(V)\subset HN_1(V)\subset ...\subset HN_k(V),$$
such that the graded quotients $gr^{HN}_i=HN_i(V)/HN_{i-1}(V)$ for
$i=1,...,k$ are semistable vector bundles and
$$\mu(gr^{HN}_1)>\mu(gr^{HN}_2)>...>\mu(gr^{HN}_k).$$
The Harder-Narasimhan filtration is unique.

A Jordan-H\"{o}lder filtration for a semistable vector bundle $V$ is
a filtration:
$$0=V_0\subset V_1\subset ...\subset V_k=V$$
such that the graded quotients $\mathrm{gr}^{V}_i=V_i/V_{i-1}$ are stable of
the same slope.

A Jordan-H\"{o}lder filtration always exists, but it is not unique
in general. The graded object $\mathrm{gr}^{V}_i=\bigoplus \mathrm{gr}^{V}_i$ does not
depend on the choice of the Jordan-H\"{o}lder filtration.

For a vector bundle $V$, define $\mu_i(V)=\mu(\mathrm{gr}^{HN}_j)$ if
$$\mathrm{rk}(HN_{j-1}(V))<i\leq \mathrm{rk}(HN_j(V)).$$ 
Obviously $\mu_1(V)\geq ...\geq \mu_k(V)$.
\begin{Lemma}\label{sub} Let $W$ be a locally free subsheaf of vector bundle
$V$ and $1\leq i\leq \mathrm{rk}(W)$. Then $\mu_i(W)\leq\mu_i(V)$. 
\end{Lemma}
\begin{proof}
If $\mu_i(W)>\mu_i(V)$, then let
$\mu_i(W)=\mu(\mathrm{gr}^{HN(W)}_j)$ and $\mu_i(V)=\mu(\mathrm{gr}^{HN(V)}_k)$. By
\cite[Lemma 1.3.3]{HL}, the canonical morphism
$HN_j(W)\hookrightarrow V\rightarrow V/HN_{k-1}(V)$ is zero, namely
$HN_j(W)\hookrightarrow HN_{k-1}(V)$, which contradicts $\mathrm{rk}(HN_j(W))\geq i>\mathrm{rk}(HN_{k-1}(V))$.
\end{proof}

Let $\mathrm{grad}(HN(V))$ denote the direct sum of the graded quotients of
the Harder-Narasimhan filtration: $\mathrm{grad}(HN(V))=\oplus \mathrm{gr}^{HN(V)}_i$.
\begin{Lemma}\label{directsum}Given vector bundles $V_1,...,V_n$, we have
$$\mathrm{grad}(HN(V_1\oplus...\oplus V_n))=\mathrm{grad}(HN(V_1))\oplus...\oplus \mathrm{grad}(HN(V_n))$$
and $\mu_i(V_1\oplus...\oplus
V_n)=\mu_i(\mathrm{grad}(HN(V_1))\oplus...\oplus \mathrm{grad}(HN(V_n)))$ for any $i$.
\end{Lemma}
\begin{proof}By induction, we only need to show the case $n=2$. Let
$$0=HN_0(V_1)\subset HN_1(V_1)\subset ...\subset HN_{k_1}(V_1)$$
$$0=HN_0(V_2)\subset HN_1(V_2)\subset ...\subset HN_{k_2}(V_2)$$
be the Harder-Narasimhan filtrations of $V_1,V_2$, respectively.

Set $0=HN_0(V_1\oplus V_2)=HN_0(V_1)\oplus HN_0(V_2)$. Assume we
have set $HN_i(V_1\oplus V_2)=HN_{i_1}(V_1)\oplus HN_{i_2}(V_2)$. We
will get $HN_{i+1}(V_1\oplus V_2)$ by the following rule:
\begin{itemize}
\item If
$$\mu(HN_{i_1+1}(V_1)/HN_{i_1}(V_1))>\mu(HN_{i_2+1}(V_2)/HN_{i_2}(V_2)),$$\\
then let $$HN_{i+1}(V_1\oplus V_2)=HN_{i_1+1}(V_1)\oplus
HN_{i_2}(V_2).$$

\item If $$\mu(HN_{i_1+1}(V_1)/HN_{i_1}(V_1))=\mu(HN_{i_2+1}(V_2)/HN_{i_2}(V_2)),$$
\\then let  $$HN_{i+1}(V_1\oplus V_2)=HN_{i_1+1}(V_1)\oplus
HN_{i_2+1}(V_2).$$

\item If $$\mu(HN_{i_1+1}(V_1)/HN_{i_1}(V_1))<\mu(HN_{i_2+1}(V_2)/HN_{i_2}(V_2)),$$
\\then let $$HN_{i+1}(V_1\oplus V_2)=HN_{i_1}(V_1)\oplus
HN_{i_2+1}(V_2).$$
\end{itemize}
 It is easy to check that the vector bundle $\mathrm{gr}^{HN(V_1\oplus V_2)}_i=HN_{i+1}(V_1\oplus V_2)/HN_i(V_1\oplus V_2)$ is semistable of
 slope $$max\{\mu(\mathrm{gr}^{HN(V_1)}_{i_1+1}),\mu(\mathrm{gr}^{HN(V_2)}_{i_2+1})\},$$
 and the slope is strictly decreasing in $i$.
 We have thus constructed the
Harder-Narasimhan filtration of $V_1\oplus V_2$. From the
construction, we also have
$$\mathrm{grad}(HN(V_1\oplus V_2))=\mathrm{grad}(HN(V_1))\oplus \mathrm{grad}(HN(V_2))$$
and $\mu_i(V_1\oplus V_2)=\mu_i(\mathrm{grad}(HN(V_1))\oplus \mathrm{grad}(HN(V_2)))$
for any $i$.
\end{proof}

\section{Weierstrass filtrations for one section}
In this section, we will consider a basic example of Weierstrass
semigroup filtrations and Weierstrass exponents.

For a Teichm\"{u}ller curve in
$\overline{\Omega\mathcal{M}}_g(2g-2)$, the universal family
$f:S\rightarrow C$ has a relative canonical bundle formula (formula
\eqref{canonical}):
$$\omega_{S/C}=f^*\mathcal{L}\otimes \mathcal{O}((2g-2)D),$$
where $D$ is the section of the unique zero. By the projection
formula (\cite[p. 124]{Ha}), 
$$f_*\omega_{S/C}=\mathcal{L}\otimes f_*\mathcal{O}((2g-2)D).$$

\subsection{Nonvarying Weierstrass semigroups in $\overline{\Omega\mathcal{M}}_g(2g-2)$}
We will consider the \emph{variation} of the Weierstrass semigroup
of
$$p=D|_F$$
along the varying fiber $F$ of $f$.
\begin{Definition}For distinct points $p_1,...,p_k$ in a Riemann surface, we define the 
 \emph{Weierstrass semigroup} $H_{p_1,\cdots,p_k}$ as follows
\begin{equation}\label{semigroup}
H_{p_1,...,p_k}:=\{(n_1,...,n_k)|h^0(n_1p_1+...+n_kp_k)=h^0(n_1p_1+...+n_kp_k-p_j)+1,\forall
j\},
\end{equation}
and we define the \emph{Weierstrass gap sequence} by the formula
$$G_{p_1,...,p_k}:=\mathbb{N}^k-H_{p_1,...,p_k}.$$
\end{Definition}

More information about the Weierstrass semigroup of one point $p$
can be found in \cite[p. 41]{ACGH85}. Using the Riemann-Roch
Theorem, we can compute the cardinality  of $G_p$, which is equal to
$g$. In fact, 
$$G_p=\{n\in\mathbb{N}:\mathrm{there\ exists\ }\omega\in H^0(C,K)\mathrm{\ with\ }\mu_p(\omega)=n-1\}.$$
This expression is closely related to the filtration we will
construct.

 We define the weight $w(H_p)$ of the Weierstrass semigroup $H_p$ to be:
\begin{equation}\label{weight}
w(H_p)=\underset{n\in G_p}{\sum}n-g(g+1)/2,
\end{equation}
which satisfies the inequality
\begin{equation}\label{weightin}
w(H_p)\leq g(g-1)/2.
\end{equation}
Moreover the equality holds if and only if $2\in H_p$.

In fact by \cite{Bu}, generic points in
$\overline{\Omega\mathcal{M}}^{odd}_g(2g-2)$ have Weierstrass gaps
$G_p=\{1, 2, 3,...,g-2,g-1,2g-1\}$, and generic points in
$\overline{\Omega\mathcal{M}}^{even}_g(2g-2)$ have Weierstrass gaps
$G_p=\{1, 2, 3,...,g-2,g,2g-1\}$.
\begin{Definition}For a Teichm\"{u}ller curve $C$ in $\overline{\Omega\mathcal{M}}_g(2g-2)$,
denote by $f:S\rightarrow C$ the universal family on $C$ and $D$ the
section of the unique zero. We say that the Weierstrass semigroup
$H_p$ (or the Weierstrass gap $G_p$ ) is \emph{nonvarying} if
$H_{D|_{F_2}}$ equals $H_{D|_{F_1}}$ (or $G_{D|_{F_2}}$ equals
$G_{D|_{F_1}}$) for any two fibers $F_1$ and $F_2$ of $f$ (including
degenerate fibers in the boundary of the moduli space).
\end{Definition}

\begin{Proposition}\label{one}For a Teichm\"{u}ller curve in
$\overline{\Omega\mathcal{M}}^{hyp}_g(2g-2)$, the Weierstrass gap
$G_p$ is nonvarying and equals $\{1, 3, 5,...,2g-5,2g-3,2g-1\}$.

If, moreover, $g\leq 5$, then
\begin{enumerate}
\item[1.] For a Teichm\"{u}ller curve in
$\overline{\Omega\mathcal{M}}^{odd}_g(2g-2)$, the Weierstrass gap
$G_p$ is nonvarying and equals $\{1, 2, 3,...,g-2,g-1,2g-1\}$.
\item[2.] For a Teichm\"{u}ller curve in
$\overline{\Omega\mathcal{M}}^{even}_g(2g-2)$, the Weierstrass gap
$G_p$ is nonvarying and equals $\{1, 2, 3,...,g-2,g,2g-1\}$.
\end{enumerate}
\end{Proposition}

\begin{proof}
The dimension $h^0(2p)$ of the fibres of $f_*\mathcal{O}(2D)$ is
an upper semicontinuous function, and $h^0(2D|_F)=2$ at smooth
fibres.
 It is also known that $$2=\mathrm{deg}(2p)\geq h^0(2p)\geq2$$ for singular fibres and $H_p$ is a
 semigroup, so $$\{2,4,...,2g-4,2g-2,2g-1,...\}\in H_p$$. Because
 $|G_p|=g$,
 and hence we get $$G_p=\mathbb{N}-H_p=\{1, 2, 3,...,g-2,g,2g-1\}.$$

For   (1) and (2), by Theorem \ref{singular}, the Teichm\"{u}ller
curve is irreducible and non-hyperelliptic. It does not have
separating nodes and $(g-1)\leq 4$, so we can use Clifford Theorem
\ref{clifford} in each of the two cases.
\begin{enumerate} 
\item[(1)]In $\overline{\Omega\mathcal{M}}^{odd}_g(2g-2)$ we have
$h^0((g-1)p)\leq 1+(g-1)/2$. If $h^0((g-1)p)=3$, the equality
implies that it is a hyperelliptic curve, hence leading to a
contradiction. So we obtain that $h^0((g-1)p)=1$.
\item[(2)] In $\overline{\Omega\mathcal{M}}^{even}_g(2g-2)$, the theta
characteristic is even, hence we have $h^0((g-1)p)=2$.
Non-hyperellipticity means that $h^0((g-2)p)=1$ for $g\leq 4$. It is
also true that $h^0(3p)=1$ for $g=5$ by Theorem \ref{singular}.
Using the Riemann-Roch Theorem we get $h^0((g+2)p)=h^0((g-2)p)+1=2$.
\end{enumerate}
\end{proof}
\subsection{Weierstrass filtrations in $\overline{\Omega\mathcal{M}}_g(2g-2)$}
On the surface $S$ we have the exact sequence
$$0\rightarrow\mathcal{O}((d-1) D)\rightarrow\mathcal{O}(d D)\rightarrow\mathcal{O}_D(d D)\rightarrow0.$$
Applying $f_*$, and using the fact that $f$ induces an isomorphism
between $D$ and $C$ ($D$ is a section), we have
$$f_*\mathcal{O}_D(d D)=\mathcal{O}_D(d D),$$
and a long exact sequncep
\begin{equation}\label{becq}
0\rightarrow f_*\mathcal{O}((d-1) D)\rightarrow f_*\mathcal{O}(d
D)\rightarrow\mathcal{O}_D(d D)\overset{\delta}{\rightarrow} R^1f_*\mathcal{O}((d-1) D)\rightarrow R^1f_*\mathcal{O}(d D).
\end{equation}
\begin{Lemma}\label{vb}$f_*\mathcal{O}(dD)$ is a vector bundle of $\mathrm{rank}(h^0(dp))$, where $p$ is $D|_F$ for a general fiber $F$. If the
Weierstrass semigroup $H_p$ is nonvarying, then
$R^1f_*\mathcal{O}(dD)$ is also a vector bundle.
\end{Lemma}
\begin{proof}By the Riemann-Roch Theorem, if $d\geq 2g$ then $h(d):=h^0(dp)=h^0(dD|_F)$ is constant for any fiber $F$.
So we get a vector bundle $f_*\mathcal{O}(dD)$ by Grauert
Semicontinuity Theorem \cite[p. 288, Corollary 12.9]{Ha}. For $0\leq
d<2g$, $f_*\mathcal{O}(dD)$ is a subsheaf of the locally free sheaf
(i.e. vector bundle) $f_*\mathcal{O}(2gD)$ by exact sequence
\eqref{becq}. Because any subsheaf of a locally free sheaf is
locally free, $f_*\mathcal{O}(dD)$ is always a vector bundle for
$d\geq 0$. By Semicontinuity Theorem \cite[p. 288, Theorem
12.7]{Ha}, the rank of $f_*\mathcal{O}(dD)$ is $h^0(dp)$.

If the Weierstrass semigroup $H_p$ is nonvarying, then the
following dimension
$h^1(dp)=h^1(dD|_F)=h^0(dD|_F)-d-\chi(\mathcal{O}_F)$ is constant
for any fiber $F$, so we get vector bundles
$f_*\mathcal{O}(dD),R^1f_*\mathcal{O}(dD)$ by Grauert semicontinuity
Theorem.
\end{proof}
If the Weierstrass semigroup is varying, in general
$R^1f_*\mathcal{O}(dD)$ is not a vector bundle\footnote{D. Chen
clarified this confusion in our first draft.}.

Let $h(d):=h^0(dp)$, for $p=D|_F$ for a general fiber $F$. Define
$$V_{h(d)}:=f_*\mathcal{O}(dD)\subset f_*\mathcal{O}((2g-2)D).$$

\begin{Remark}The definition is reasonable:
if $h^0(dD|_F)=h^0((d+1)D|_F)$, then
$f_*\mathcal{O}(dD)=f_*\mathcal{O}((d+1)D)$.
\end{Remark}
Thus we get a filtration of the vector bundle
$f_*\mathcal{O}((2g-2)D)$:
$$0\subset V_1\subset V_2 \subset... \subset V_g=f_*\mathcal{O}((2g-2)D).$$
It is an example of \emph{Weierstrass filtrations}. If the
Weierstrass semigroup $H_p$ is nonvarying, It is also an example of
\emph{Weierstrass semigroup filtrations}.

Denote by $d_i$ the $i$-th element in $H_p$, for $p=D|_F$ for a
general fiber $F$.

\begin{Lemma}\label{degree}The graded quotient $V_i/V_{i-1}$ is a line bundle of degree at most $\frac{-d_i}{2g-1}\mathrm{deg} \mathcal{L}$. If the
Weierstrass semigroup $H_p$ is nonvarying, then the degree of
$V_i/V_{i-1}$ equals $\frac{-d_i}{2g-1}\mathrm{deg} \mathcal{L}$.
\end{Lemma}

\begin{proof}
For $d_i\in H_p$ the $i$-th element in $H_p$,
$\mathrm{rk}f_*\mathcal{O}((d_i-1) D)=\mathrm{rk}f_*\mathcal{O}(d_i
D)-1$. Therefore
$$V_i=f_*\mathcal{O}(d_i D),\quad V_{i-1}=f_*\mathcal{O}((d_i-1) D).$$
We have a long exact sequence by \eqref{becq}:
$$0\rightarrow f_*\mathcal{O}((d_i-1) D)\rightarrow f_*\mathcal{O}(d_i D)\rightarrow\mathcal{O}_D(d_i D)\overset{\delta}{\rightarrow}R^1f_*\mathcal{O}((d-1)D).$$
$\mathrm{ker}(\delta)$ is a subsheaf of the line bundle
$\mathcal{O}_D(d_i D)$, so it is also a line bundle whose degree
at most $\mathrm{deg}(\mathcal{O}_D(d_i D))$. Then by
the formula \eqref{intersection}:
\begin{align*}
    \mathrm{deg}(V_i/V_{i-1}) & =\mathrm{deg}(f_*\mathcal{O}(d_i D))-\mathrm{deg}(f_*\mathcal{O}((d_i-1)
D))\\
     & \leq \mathrm{deg}(\mathcal{O}_D(d_i
D))=d_iD^2=\frac{-d_i}{2g-1}\mathrm{deg} \mathcal{L},
\end{align*}p
Since subsheaves of a locally free sheaf are locally free, we deduce that
$\mathrm{ker}(\delta)$ and $\mathrm{im}(\delta)$ are both locally
free.

If, moreover, the Weierstrass semigroup $H_p$ is nonvarying, then
$R^1f_*\mathcal{O}((d-1)D)$ is locally free by Lemma \ref{vb}.
$\mathrm{ker}(\delta)$ and $\mathrm{im}(\delta)$ are both locally
free. $\mathcal{O}_D(d_i D)$ is a line bundle, so
$\mathrm{im}(\delta)=\mathcal{O}_D(d_i D)/\mathrm{ker}(\delta)$ is
zero. We have a short exact sequence
$$0\rightarrow f_*\mathcal{O}((d_i-1) D)\rightarrow f_*\mathcal{O}(d_i D)\rightarrow\mathcal{O}_D(d_i D)\overset{\delta}\rightarrow0,$$
and the degree of $V_i/V_{i-1}$ equals
$\frac{-d_i}{2g-1}\mathrm{deg} \mathcal{L}$.
\end{proof}
We get a filtration of $f_*\omega_{S/C}=\mathcal{L} \otimes
f_*\mathcal{O}((2g-2)D)$:
$$0\subset \mathcal{L}\otimes V_1\subset \mathcal{L}\otimes V_2\subset... \subset \mathcal{L}\otimes V_g=\mathcal{L} \otimes f_*\mathcal{O}((2g-2)D)=f_*\omega_{S/C}.$$
\begin{Definition}If the Weierstrass semigroup  $H_p$ is nonvarying, we define the
$i$-th \emph{Weierstrass exponent} $w_i$ as follow:
$$w_i=\mathrm{deg}(\mathcal{L}\otimes V_i/\mathcal{L}\otimes V_{i-1})/\mathrm{deg}(\mathcal{L})=1-\frac{d_i}{2g-1}=\frac{2g-1-d_i}{2g-1}.$$
\end{Definition}
\begin{Remark}The sum of Weierstrass exponents is
 $\mathrm{deg}(f_*\omega_{S/C})/\mathrm{deg}(\mathcal{L})$, which
equals the sum of Lyapunov exponents by Theorem \ref{sumly}.
\end{Remark}

\begin{Theorem}\label{2g-2}Let $k_i$ be the $i$-th
largest number in the Weierstrass gap $G_p$. For a Teichm\"{u}ller curve in
$\overline{\Omega\mathcal{M}}^{hyp}_g(2g-2)$, the Weierstrass
exponent $w_i$ is $\frac{k_i}{2g-1}$ and the Weierstrass gap $G_p$
equals $\{1, 3, 5,...,2g-5,2g-3,2g-1\}$.

If, moreover, $g\leq 5$, then:
\begin{enumerate}
\item[1.]For a Teichm\"{u}ller curve in
$\overline{\Omega\mathcal{M}}^{odd}_g(2g-2)$, the Weierstrass
exponent $w_i$ is $\frac{k_i}{2g-1}$ and the Weierstrass gap $G_p$
equals $\{1, 2, 3,...,g-2,g-1,2g-1\}$, and $f_*{\omega_{S/C}}$
splits into direct sum of line bundles.
\item[2.] For a Teichm\"{u}ller curve in
$\overline{\Omega\mathcal{M}}^{even}_g(2g-2)$, the Weierstrass
exponent $w_i$ is $\frac{k_i}{2g-1}$ and the Weierstrass gap $G_p$
equals $\{1, 2, 3,...,g-2,g,2g-1\}$.
\end{enumerate}
\end{Theorem}
\begin{proof}Proposition \ref{one} tells us that these Teichm\"{u}ller
curves have nonvarying Weierstrass semigroups, and for $d_i\in
H_p$, $2g-1-d_i\in G_p$, hence we obtain the result by applying
Lemma \ref{degree} and Theorem \ref{odd}.
\end{proof}
The weight formula \eqref{weight} gives the sum formula
\begin{equation}\label{onesum}
\overset{g}{\underset{i=1}{\sum}} w_i=\frac{1}{2g-1}(\underset{n\in
G_p}{\sum} n)=\frac{1}{2g-1}{w(H_p)+\frac{g(g+1)}{2(2g-1)}}.
\end{equation}
It has maximal value $\frac{g^2}{(2g-1)}$ by the inequality
\eqref{weightin}, where the equality holds if and only if $2\in
H_p$.
\begin{Remark}
If the Weierstrass semigroup $H_p$ is varying, the sum formula
\eqref{onesum} also gives a upper bound of $\mathrm{deg}
f_*\omega_{S/C}/\mathrm{deg}(\mathcal{L})$ by using Lemma
\ref{degree}.
\end{Remark}
\begin{Corollary}
If the Weierstrass semigroup $H_p$ is nonvarying, then
$$0\subset V_1\subset V_2 \subset... \subset V_g=f_*\mathcal{O}((2g-2)D)$$
is the Harder-Narasimhan filtration.
\end{Corollary}
\begin{proof}Lemma \ref{degree} tells us that $\mathrm{deg}(V_i/V_{i-1})>\mathrm{deg}(V_{i+1}/V_i)$, hence we conclude the result by the uniqueness
of the Harder-Narasimhan filtration.
\end{proof}

\subsection{$f_*\mathcal{O}_{aD}(dD)$ and Splitting Lemma I}
We can get more information by analyzing the exact sequence
\begin{equation}\label{basic}
0\rightarrow f_*\mathcal{O}((d-a) D)\rightarrow f_*\mathcal{O}(d
D)\rightarrow f_*\mathcal{O}_{aD}(d D)\overset{\delta}{\rightarrow}R^1f_*\mathcal{O}((d-a) D)\rightarrow.
\end{equation}
Here $\mathrm{ker}(\delta)$ is controlled by
$f_*\mathcal{O}_{aD}(dD)$ ($=\mathcal{O}_{aD}(dD)$ via the equation
\eqref{direct}), which has a good filtration.

\begin{Lemma}\label{hn}The Harder-Narasimhan filtration of
$f_*\mathcal{O}_{aD}(dD)$ is
$$0\subset f_*\mathcal{O}_{D}((d-a+1)D)\subset...\subset f_*\mathcal{O}_{(a-1)D}((d-1)D)\subset f_*\mathcal{O}_{aD}(dD).$$p
and the direct sum of the graded quotients of this filtration is
$$\mathrm{grad}(HN(f_*\mathcal{O}_{aD}(dD)))=\bigoplus_{i=0}^{a-1}
\mathcal{O}_{D}((d-i)D).$$
\end{Lemma}
\begin{proof}
Because $f$ is an isomorphism between $D$ and $C$ ($D$ is a
section),
$$f_*\mathcal{O}_{D}(j D)=\mathcal{O}_{D}(j
D),\quad  R^1f_*\mathcal{O}_{D}(j D)=0.$$ From the exact sequence
$$0\rightarrow \mathcal{O}_{(i-1)D}((j-1)D) \rightarrow \mathcal{O}_{iD}(jD) \rightarrow \mathcal{O}_{D}(jD) \rightarrow 0$$
with $1 \leq i\leq a,\ d-a+1 \leq i\leq d$, we get the long exact
sequence
$$0\rightarrow f_*\mathcal{O}_{(i-1)D}((j-1)D) \rightarrow f_*\mathcal{O}_{iD}(jD) \rightarrow \mathcal{O}_{D}(jD) $$
$$\qquad\qquad\quad   \rightarrow R^1f_*\mathcal{O}_{(i-1)D}((j-1)D) \rightarrow R^1f_*\mathcal{O}_{iD}(jD) \rightarrow R^1f_*\mathcal{O}_{D}(jD).$$
By induction, we have
\begin{equation}\label{direct}
f_*\mathcal{O}_{iD}(j D)=\mathcal{O}_{iD}(j D),\quad
R^1f_*\mathcal{O}_{iD}(j D)=0,
\end{equation}
and the exact sequence
$$0\rightarrow f_*\mathcal{O}_{(i-1)D}((j-1)D) \rightarrow f_*\mathcal{O}_{iD}(jD) \rightarrow \mathcal{O}_{D}(jD) \rightarrow 0.$$

 Because $D^2<0$, we obtain a filtration with a property that the
degree of graded quotient line bundles are strictly decreasing. By
the uniqueness of the Harder-Narasimhan filtration, we conclude the
desired result.
\end{proof}
The filtration can be used to describe the structure of special
quotients.
\begin{Lemma}[Splitting Lemma I]\label{I} If the dimensions $$h^0(dp),h^0((d-m)p),h^0((d-m-n+1)p)$$ are constant for any fiber,
and satisfy $$h^0(dp)=h^0(d-m)p)+m=h^0((d-m-n+1)p)+m,$$ then
$f_*\mathcal{O}_{aD}$ splits in $f_*\mathcal{O}_{(a+b)D},\ a\leq m,\
b<n$. In particular, for $n\geq 2$,
$$f_*\mathcal{O}(dD)/f_*\mathcal{O}((d-m)D)=\bigoplus_{i=0}^{m-1} \mathcal{O}_{D}((d-i)D).$$
\end{Lemma}
\begin{proof}We have the following commutative diagram
$$
\xymatrix{
f_*\mathcal{O}((d-m)D)\ar@{^{(}->}[r]&f_*\mathcal{O}((d-m+a)D) \ar@{->>}[r]^{\vartheta}&  f_*\mathcal{O}_{aD}((d-m+a)D)  \ar@/_10mm/[d]^{\phi}\\
f_*\mathcal{O}((d-m-b)D)\ar@{^{(}->}[r]\ar@{->}[u]^{\theta||}&f_*\mathcal{O}((d-m+a)D)\ar@{->}[u]^{||} \ar@{->}[r]^{\psi}   &  f_*\mathcal{O}_{(a+b)D}((d-m+a)D)   \ar@{->>}[u]^{\varphi}\\
            &  & \ar@{^{(}->}[u]f_*\mathcal{O}_{bD}((d-m)D) }
$$
Because $h^0(dp)=h^0(d-m)p)+m$, a similar argument as in Lemma
\ref{degree}  implies that $\vartheta$ is surjective. Moreover
$h^0((d-m)p)=h^0((d-m-n+1)p)$ implies that $\theta$ is an
isomorphism by Corollary \ref{iso}.

Thus the image of $\psi$ is the same as the image of $\vartheta$,
that is $f_*\mathcal{O}_{aD}((d-m+a)D)$. So there is a $\phi$ with
$\varphi\phi=id$, hence $f_*\mathcal{O}_{aD}((d-m+a)D)$ splits in
$f_*\mathcal{O}_{(a+b)D}((d-m+a)D),\ a\leq m,\ b<n$.

If $n\geq 2$, $f_*\mathcal{O}_{aD}$ splits in
$f_*\mathcal{O}_{(a+1)D}$ for $a\leq m$. By induction,
$f_*\mathcal{O}_{mD}$ splits into a direct sum of line bundles.

 Thus by the uniqueness of the Harder-Narasimhan filtration and Lemma \ref{hn}, we have
$$f_*\mathcal{O}(dD)/f_*\mathcal{O}((d-m)D)=f_*\mathcal{O}_{mD}(dD)=\bigoplus^{m-1}_{i=0} \mathcal{O}_{D}((d-i)D).$$
\end{proof}
\begin{Theorem}\label{odd}
If $g\leq 5$, the sheaf $f_*\omega_{S/C}$ of a Teichm\"{u}ller
curve $C$ in $\overline{\Omega\mathcal{M}}^{odd}_g(2g-2)$ splits
into a direct sum of line bundles.
\end{Theorem}
\begin{proof}Theorem \ref{one} gives
$h^0(2(g-1)p)=h^0((g-1)p)+g-1=h^0((g-2)p)+g-1$. We have
\begin{align*}
 f_*\omega_{S/C} &= \mathcal{L}\oplus(\mathcal{L} \otimes (f_*\mathcal{O}((2g-2)D)/f_*\mathcal{O}))\\
 &= \mathcal{L}\oplus(\mathcal{L} \otimes
 (f_*\mathcal{O}((2g-2)D)/f_*\mathcal{O}((g-1)D)))\\
&=\mathcal{L}\oplus(\bigoplus^{g-1}_{i=0}
\mathcal{O}_{D}((2g-2-i)D)\otimes\mathcal{L}),
\end{align*}
where the first equality is by Equation \eqref{VHS} and the last
equality follows from Lemma \ref{I}.
\end{proof}
\section{Weierstrass filtrations for several sections}
In this section we will define three kinds of filtrations:
Weierstrass filtrations, Weierstrass semigroup filtrations  and
Weierstrass pair filtrations. The first one together with the upper
bound Lemma \ref{up} is used to get coarse information about the
upper bound of the sum. The second and the third can be used to get
more precise information about each quotient.

\subsection{Weierstrass filtrations}
 From the exact sequence
$$0\rightarrow f_*\mathcal{O}(d_1D_1+...+d_kD_k)\rightarrow f_*\mathcal{O}(m_1D_1+...+m_kD_k)=f_*(\omega_{S/C})\otimes\mathcal{L}^{-1}$$
and the fact that all subsheaves of a locally free sheaf on a curve
are locally free, we deduce that $f_*\mathcal{O}(d_1D_1+...+d_kD_k)$
is a vector bundle. By Semicontinuity Theorem \cite[p. 288, Theorem
12.7]{Ha}, its rank is $h^0(d_1p_1+...+d_kp_k)$, where $p_i=D_i|_F$
and $F$ is a generic fiber.

For $1\leq i\leq g$, we define
$$W_i =\{(d_1,...,d_k)|h^0(d_1p_1+...+d_kp_k)=i\mathrm{\ for\ a \ general\ fiber}\}$$
and nonvarying sets
$$WS_i =\{(d_1,...,d_k)|h^0(d_1p_1+...+d_kp_k)=i \mathrm{\ for\  all\ fibers\ (incl.\ boundary\ points)}\}.$$

For each element $(d_1,...,d_k)$ in $W_i$ or $WS_i$,
$f_*\mathcal{O}(d_1D_i+...+d_kD_k)$ is a rank $i$ vector bundle. If
$(d_1,...,d_k)\in WS_i$ then $R^1f_*\mathcal{O}(d_1D_i+...+d_kD_k)$
is also a vector bundle by Grauert Semicontinuity Theorem \cite[p.
288, Corollary 12.9]{Ha}.

We also define $(d_1,...,d_k)\leq (d'_1,...,d'_k)$ if $d_i\leq d'_i$
for all $i$ and define $(d_1,...,d_k)<(d'_1,...,d'_k)$ if $d_i\leq
d'_i$ but not all $d_i=d'_i$.

Denote by $\underline{d}_i$ the tuple $(d_{i1},...,d_{ik})$.
\begin{Definition} We define the set of \emph{Weierstrass filtrations}  as follows:
$$WF=\{\{\underline{d}_i\}|\mathrm{\ at\ most\ one\ }\underline{d}_i \in
W_i \mathrm{\ for\ } 1\leq i\leq g.\mathrm{\ If\ }i<j,\mathrm{\
then\ }\underline{d}_i<\underline{d}_j\}.$$
\end{Definition}
An element $\{\underline{d}_i\}\in WF$ is a filtration of vector
bundles of $f_*\mathcal{O}(m_1D_1+...+m_kD_k)$.
$$0\subset...\subset f_*\mathcal{O}(d_{i1}D_1+...+d_{ik}D_k)\subset... \subset f_*\mathcal{O}(m_1D_1+...+m_kD_k).$$

As what we have shown in Lemma \ref{degree}, the filtration dose not
give us the desired properties to compute the degree of graded
quotient vector bundles because in general
$R^1f_*\mathcal{O}(d_1D_i+...+d_kD_k)$ is not locally free. So in
many cases, we need the nonvarying assumption to get more
information.

\begin{Definition} We define the set of \emph{Weierstrass semigroup filtrations} as follows:
$$WSF=\{\{\underline{d}_i\}|\mathrm{for\ }1\leq i\leq g,\mathrm{\ there\ is\ only\ one\ } \underline{d}_i \in
WS_i\cap H_{p_1,...,p_k}, \mathrm{\ and\ }
\underline{d}_i<\underline{d}_{i+1}\},$$ where $H_{p_1,...,p_k}$ is
the Weierstrass semigroup for a general fiber.
\end{Definition}

Assuming that the Weierstrass semigroup is nonvarying, we will use
the Weierstrass semigroup filtration to define Weierstrass exponents
in the next section.
\begin{Example}For a
Teichm\"{u}ller curve in the hyperelliptic locus of the stratum
$\overline{\Omega\mathcal{M}}_5(2,2,4)$, and denote by $(4,3)$ the
orders of singularities of the corresponding quadratic
differentials.  The Weierstrass semigroup filtration we will
construct in Proposition \ref{hyperfiltration} for defining
Weierstrass exponents is
$$\{(0,0,0),(0,0,2),(0,0,4),(1,1,4),(2,2,4)\},$$
that is
$$0\subset f_*\mathcal{O}\subset  f_*\mathcal{O}(2D_3)\subset
f_*\mathcal{O}(4D_3)\subset f_*\mathcal{O}(D_1+D_2+4D_3)\subset
f_*\mathcal{O}(2D_1+2D_2+4D_3).
$$
Moreover, there is another Weierstrass semigroup filtration
$$\{(0,0,0),(1,1,0),(2,2,0),(2,2,2),(2,2,4)\}\},$$
that is
$$0\subset f_*\mathcal{O}\subset f_*\mathcal{O}(D_1+D_2)\subset
f_*\mathcal{O}(2D_1+2D_2)\subset
f_*\mathcal{O}\Big(\overset{3}{\underset{i=1}{\sum}}2D_i\Big)\subset
f_*\mathcal{O}(2D_1+2D_2+4D_3).
$$
\end{Example}

Under some weak assumptions, the following filtration is also useful
for  computational and theoretical reasons.
\begin{Definition}  We define  the set of \emph{Weierstrass pair filtrations} as follows:
$$WPF=\{\{\{\underline{d}_i,\underline{d}'_{i+1}\}\}|\mathrm{for\ }1\leq i\leq g-1,\ \underline{d}_i \in WS_i,\ \underline{d}'_i \in W_i,$$
$$\underline{d}'_i<\underline{d}_i,\ \underline{d}'_{i+1}=\underline{d}_i+(0,...,1,...,0)\}.$$
\end{Definition}
To define these we need to verify that the exact sequence
$$0\rightarrow f_*\mathcal{O}(\underline{d}_i)\rightarrow f_*\mathcal{O}(\underline{d}'_{i+1})\rightarrow\mathcal{O}_{\underline{d}'_{i+1}-\underline{d}_i}(\underline{d}'_{i+1})\overset{\delta}{\rightarrow}R^1f_*\mathcal{O}(\underline{d}_i)\rightarrow R^1f_*\mathcal{O}(\underline{d}'_{i+1})$$
satisfies $\mathrm{rk}f_*\mathcal{O}(\underline{d}_i)+1= \mathrm{rk}
f_*\mathcal{O}(\underline{d}'_{i+1})$,
$\mathcal{O}_{\underline{d}'_{i+1}-\underline{d}_i}(\underline{d}'_{i+1})$
is a line bundle and that $R^1f_*\mathcal{O}(\underline{d}_i)$ is
locally free. Note also that $\underline{d}'_i<\underline{d}_i$
implies
$f_*\mathcal{O}(\underline{d}'_i)=f_*\mathcal{O}(\underline{d}_i)$.
Thus an element $\{\{\underline{d}_i,\underline{d}'_{i+1}\}\}\in
WPF$ is a filtration of vector bundles:
$$0\subset...\subset f_*\mathcal{O}(\underline{d}_i) \subset f_*\mathcal{O}(\underline{d}'_{i+1})=f_*\mathcal{O}(\underline{d}_{i+1})\subset... \subset f_*\mathcal{O}(m_1D_1+...+m_kD_k).$$

\begin{Example}There is a Weierstrass pair filtration in the proof of the stratum
$\overline{\Omega\mathcal{M}}^{even}_4(4,2)$:
$$\{\{(1,1),(2,1)\},\{(3,1),(3,2)\},\{(3,2),(4,2)\}\},$$
that is
$$0\subset f_*\mathcal{O}(D_1+D_2)\subset
f_*\mathcal{O}(2D_1+D_2)=f_*\mathcal{O}(3D_1+D_2)$$ $$\subset
f_*\mathcal{O}(3D_1+2D_2)=f_*\mathcal{O}(3D_1+2D_2)\subset
f_*\mathcal{O}(4D_1+2D_2).
$$
\end{Example}
\begin{Remark}The relationship among those sets of filtrations is
$$WSF\subset WPF\subset WF.$$
\end{Remark}
\subsection{Splitting Lemma II and Upper Bound Lemma}
The next Lemma describes a splitting structure of the quotient:
\begin{Lemma}[Splitting Lemma II]\label{II}If $h^0(\sum d_ip_i)$ and $h^0(\sum (d_i-a_i)p_i)$ are constant for any fiber, and satisfy $h^0(\sum d_ip_i)=h^0(\sum (d_i-a_i)p_i)+\sum a_i$, then
$$f_*\mathcal{O}\Big(\sum d_iD_i\Big)/f_*\mathcal{O}\Big(\sum (d_i-a_i)D_i\Big)=\bigoplus f_*\mathcal{O}_{a_iD_i}(d_iD_i).$$
\end{Lemma}
\begin{proof}From the exact sequence
$$0\rightarrow\mathcal{O}\Big(\sum (d_i-a_i)D_i\Big)\rightarrow\mathcal{O}\Big(\sum d_iD_i\Big)\rightarrow\mathcal{O}_{\sum a_iD_i}\Big(\sum d_iD_i\Big)\rightarrow 0,$$
we get the long exact sequence
$$0\rightarrow f_*\mathcal{O}\Big(\sum (d_i-a_i)D_i\Big)\rightarrow f_*\mathcal{O}\Big(\sum d_iD_i\Big)\rightarrow f_*\mathcal{O}_{\sum a_iD_i}\Big(\sum d_iD_i\Big)$$
$$\qquad\qquad\qquad\qquad\overset{\delta}{\rightarrow}R^1f_*\mathcal{O}\Big(\sum (d_i-a_i)D_i\Big)\rightarrow R^1f_*\mathcal{O}\Big(\sum d_iD_i\Big)\rightarrow 0.$$
Because $\mathrm{ker}(\delta)$ and $\mathrm{im}(\delta)$ are both
locally free and because
$$\mathrm{rk}\Big(R^1f_*\mathcal{O}\Big(\sum (d_i-a_i)D_i\Big)\Big)=\mathrm{rk}\Big(R^1f_*\mathcal{O}\Big(\sum
d_iD_i\Big)\Big),$$
we get
$$0\rightarrow f_*\mathcal{O}\Big(\sum (d_i-a_i)D_i\Big)\rightarrow f_*\mathcal{O}\Big(\sum d_iD_i\Big)\rightarrow f_*\mathcal{O}_{\sum a_iD_i}\Big(\sum d_iD_i\Big)\rightarrow 0.$$
Since $D_i\cdot D_j=0$ for $i\neq j$, we have
\begin{align*}
    f_*\mathcal{O}_{\sum a_iD_i}\Big(\sum d_iD_i\Big) &=f_*\Big(\bigoplus \mathcal{O}_{a_iD_i}\Big(\sum d_iD_i\Big)\Big)\\
     &=f_*(\bigoplus \mathcal{O}_{a_iD_i}(d_iD_i))=\bigoplus f_*\mathcal{O}_{a_iD_i}(d_iD_i).
\end{align*}
\end{proof}
We often use the above Lemma with Lemma \ref{I}.
\begin{Corollary}
For a Teichm\"{u}ller curve in
$\overline{\Omega\mathcal{M}}_g(2k_1,...,2k_m)$ with theta
characteristic $h^0(\sum k_ip_i)=1$ for every fiber, the sheaf
$f_*{\omega_{S/C}}$ splits into a direct sum of line bundles.
\end{Corollary}

\begin{proof}Because
$$h^0\Big(\sum 2k_ip_i\Big)=h^0\Big(\sum k_ip_i\Big)+g-1=h^0\Big(\sum k_ip_i\Big)+\sum k_i,$$
we can apply Lemma \ref{II}. Because
$$h^0\Big(\Big(\sum k_ip_i\Big)+p_j\Big)-1=h^0\Big(\sum k_ip_i\Big)=1=h^0\Big(\Big(\sum k_ip_i\Big)-p_j\Big).$$
it suffices to apply Lemma \ref{I} to each $D_j$.
\end{proof}

\begin{Corollary}\label{WSF} Each graded quotient of a filtration in $WSF$(resp. in $WPF$) is a line bundle
whose degree is $-\frac{d_j}{m_j+1}\mathrm{deg}(\mathcal{L})$ of for
some $d_j,m_j$. So the degree of $f_*\mathcal{O}(m_1D_1+...+m_kD_k)$
(i.e. the sum of the degree of those line bundles) can be computed.
All filtrations in $WSF$ have the same sum
$\mathrm{deg}f_*\omega_{S/C}-g \mathrm{deg}(\mathcal{L})$.
\end{Corollary}

\begin{proof}For any filtration in $WSF$, we
have the exact sequence
$$WSF:0\rightarrow f_*\mathcal{O}(\underline{d}_i-D_j)\rightarrow f_*\mathcal{O}(\underline{d}_i)\rightarrow\mathcal{O}_{D_j}(\underline{d}_i)\rightarrow 0,$$
where $\underline{d}_{i-1}\leq\underline{d}_i-D_j$. Similarly for $WPF$ we have the exact sequence
$$WPF:0\rightarrow f_*\mathcal{O}(\underline{d}_i)\rightarrow f_*\mathcal{O}(\underline{d}'_{i+1})\rightarrow\mathcal{O}_{\underline{d}'_{i+1}-\underline{d}_i}(\underline{d}'_{i+1})\rightarrow 0.$$

Therefore each graded quotient of the filtration is a line bundle
$\mathcal{O}_{D_j}(\underline{d}_i)$ (resp.\
$\mathcal{O}_{\underline{d}'_{i+1}-\underline{d}_i}(\underline{d}'_{i+1})$),
whose degree is
$d_jD_j^2=-\frac{d_j}{m_j+1}\mathrm{deg}(\mathcal{L})$ for some
$d_j$ by formula (\ref{intersection}).

The degree of $f_*\mathcal{O}(m_1D_1+...+m_kD_k)$ is the sum of the
degree of those line bundles.

For the last part of the claim, the sum is
$$\mathrm{deg} f_*\mathcal{O}(m_1D_1+...+m_kD_k)=\mathrm{deg}f_*\omega_{S/C}-g
\mathrm{deg}(\mathcal{L}).$$
\end{proof}
\begin{Example}An explicit formula about the sum has been established in Corollary \ref{EKZ}
for a Teichm\"{u}ller curve in the hyperelliptic locus.
\end{Example}
When $h^0(\sum d_ip_i)\leq h^0(\sum (d_i-a_i)p_i)+\sum a_i$,  we can
get an upper bound for the quotient by the properties of the
Harder-Narasimhan filtration, even if $h^0(\sum d_ip_i)$ and
$h^0(\sum (d_i-a_i)p_i)$ are varying.
\begin{Lemma}[Upper Bound Lemma]\label{up} Let $V=f_*\mathcal{O}(\sum d_iD_i)/f_*\mathcal{O}(\sum (d_i-a_i)D_i)$ and $r=h^0(\sum
d_ip_i)-h^0(\sum (d_i-a_i)p_i)$, where $p_i$ is $D_i|_F$ for a
general fiber $F$. We order degrees of line bundles in the set
$$\{\mathcal{O}_{D_i}((d_i-j)D_i)|1\leq i\leq k, 0\leq j\leq a_i-1\}$$
decreasing (counted with multiplicity) as $b_1\geq b_2\geq...\geq
b_{a_1+...+a_k}$. Then
$$\mathrm{deg}(V)\leq b_1+b_2+...+b_r.$$
\end{Lemma}

\begin{proof}
By Lemma \ref{hn} and Lemma \ref{directsum}, the graded sum of the
Harder-Narasimhan filtration of $\oplus_i
f_*\mathcal{O}_{a_iD_i}(d_iD_i)$ is

\begin{align*}
   \mathrm{grad}(HN(\bigoplus_i
f_*\mathcal{O}_{a_iD_i}(d_iD_i)) &=\bigoplus_i\mathrm{grad}(HN(f_*\mathcal{O}_{a_iD_i}(d_iD_i)))\\
     &=\bigoplus_i\bigoplus^{a_i-1}_{j=0}\mathcal{O}_{D_i}((d_i-j)D_i).
\end{align*}p
So
$\mu_l(\bigoplus_if_*\mathcal{O}_{a_iD_i}(d_iD_i))=\mu_l(\bigoplus_i\bigoplus^{a_i-1}_{j=0}\mathcal{O}_{D_i}((d_i-j)D_i))$
for all $l$. From the proof of Splitting Lemma \ref{II}, the kernel
$\mathrm{ker}(\delta)$ is a locally free subsheaf of rank
$r=h^0(\sum d_ip_i)-h^0(\sum (d_i-a_i)p_i)$:
$$V=f_*\mathcal{O}\Big(\sum d_iD_i\Big)/f_*\mathcal{O}\Big(\sum (d_i-a_i)D_i\Big)\subset\bigoplus_i f_*\mathcal{O}_{a_iD_i}(d_iD_i).$$
So by Lemma \ref{sub}, $$\mathrm{deg}(V)\leq
\sum^{r}_{l=1}\mu_l(\bigoplus_if_*\mathcal{O}_{a_iD_i}(d_iD_i))=
\sum^{r}_{l=1}\mu_l(\bigoplus_i\bigoplus^{a_i-1}_{j=0}\mathcal{O}_{D_i}((d_i-j)D_i)).$$
It is obvious that
$\mu_l(\oplus_i\oplus^{a_i-1}_{j=0}\mathcal{O}_{D_i}((d_i-j)D_i))$
is the $l$-th largest degree in the set of line bundles
$$\{\mathcal{O}_{D_i}((d_i-j)D_i)|1\leq i\leq k, 0\leq j\leq a_i-1\}.$$
\end{proof}
\begin{Corollary}\label{iso}If $h^0(\sum d_ip_i)=h^0(\sum (d_i-a_i)p_i)$ holds in a general fiber, then $f_*\mathcal{O}(\sum d_iD_i)=f_*\mathcal{O}(\sum (d_i-a_i)D_i)$.
\end{Corollary}

\subsection{Application to the sum of Lyapunov exponents}
The existence of Weierstrass semigroup (pair) filtrations is
convenient for computation.

\begin{Corollary}\label{sum}
If there exists a Weierstrass semigroup (pair) filtration for a
Teichm\"{u}ller curve, then we can compute the sum of Lyapunov
exponents. Moreover, the denominator of the sum of Lyapunov
exponents divides $(m_1+1)\cdots(m_k+1)$.
\end{Corollary}
\begin{proof}
The sum of Lyapunov exponents is
$L(C)=\mathrm{deg}f_*\omega_{S/C}/\mathrm{deg}(\mathcal{L})$ by
Theorem \ref{sumly}. The sum of the degree of each graded quotient
line bundles of the Weierstrass semigroup (pair) filtration is
$\mathrm{deg}f_*\omega_{S/C}-g \mathrm{deg}(\mathcal{L})$ by
Corollary \ref{WSF}.

Because $D^2_i/\mathrm{deg} \mathcal{L}=-\frac{1}{m_i+1}$, each
graded quotient has denominator $(m_i+1)$, hence we get the second
claim.
\end{proof}

In many cases, because of the absence of Weierstrass semigroup
(pair) filtrations, we can not  get more precise information about
$f_*\omega_{S/C}$. But the partial filtration is enough to give some
upper bound for it. Using a coarse filtration, we have:
\begin{Corollary} The sum of Lyapunov exponents for a Teichm\"{u}ller
curve $C$ in the stratum
$\overline{\Omega\mathcal{M}}_g(m_1,...,m_k)$ satisfied the
inequality\footnote{We have obtained a better upper bound $L(C)\leq \frac{g+1}{2}$ in
\cite{YZ}.}
$$L(C)\leq \frac{3g}{4}-\frac{1}{8}\Big(-2+\underset{m_i\mathrm{\ even}}{\sum}\frac{m_i}{m_i+1}+\underset{m_i\mathrm{\ odd}}{\sum}1\Big).$$
\end{Corollary}
\begin{proof}
There is a rank $g-1$ subbundle
$$ f_*\mathcal{O}(m_1D_1+...+m_kD_k)/\mathcal{O}_C\subset\bigoplus_i f_*\mathcal{O}_{m_iD_i}(m_iD_i).$$
It is obvious by Lemma \ref{hn} and Lemma \ref{directsum} that
$$\mathrm{grad}(HN(\underset{i}{\oplus}
f_*\mathcal{O}_{m_iD_i}(m_iD_i))=\bigoplus_i\bigoplus^{m_i}_{j=1}\mathcal{O}_{D_i}(jD_i).$$
By Lemma \ref{up}, we want to get the sum of the largest degrees of
$g-1$ line bundles in the set
$$\{\mathcal{O}_{D_i}(jD_i)|1\leq i\leq k, 0\leq j\leq m_i\}$$
The sum of the largest degrees of $n_i$ line bundles in the set
$\{\mathcal{O}_{D_i}(jD_i)|0\leq j\leq m_i\}$ is the sum of the
degrees of line bundles in the set $\{\mathcal{O}_{D_i}(jD_i)|0\leq
j\leq n_i\}$. It is
$\Big(-\frac{n_i(n_i+1)}{2(m_i+1)}\Big)\mathrm{deg}\mathcal{L}$. We
get
 $$\mathrm{deg}(f_*\mathcal{O}(m_1D_1+...+m_kD_k)/\mathcal{O}_C)\leq \underset{\sum n_i=g-1}{\mathrm{max}}\Big\{
-\sum\frac{n_i(n_i+1)}{2(m_i+1)}\Big\}\mathrm{deg}\mathcal{L}.$$

Because $\sum n_i=g-1=(\sum m_i)/2$, there must be some $n_j$ such
that $n_j\leq(m_j-1)/2$ if there is $n_i>(m_i+1)/2$. In this case,
we have
$$\mathrm{deg}(\mathcal{O}_{D_i}(n_iD_i))< -\frac{1}{2}\mathrm{deg}\mathcal{L}\leq
\mathrm{deg}(\mathcal{O}_{D_j}((n_j+1)D_j)).$$ There must be some
$n_i$ such that $n_i\geq(m_i+1)/2$ if there is $n_j<(m_j-1)/2$, and
in this case, we have
$$\mathrm{deg}(\mathcal{O}_{D_i}(n_iD_i))\leq -\frac{1}{2}\mathrm{deg}\mathcal{L}<
\mathrm{deg}(\mathcal{O}_{D_j}((n_j+1)D_j)).$$

We can increase the value $-\sum\frac{n_i(n_i+1)}{2(m_i+1)}$ by
changing $n_i$ to $n_i-1$ and $n_j$ to $n_j+1$ in both cases.

So when $-\sum\frac{n_i(n_i+1)}{2(m_i+1)}$ reaches the maximum, we
know that $n_i=m_i/2$ when $m_i$ is even, and $n_i=(m_i-1)/2$ or
$(m_i+1)/2$ when $m_i$ is odd, with the property that
$$\mathrm{Card}(\{n_i=(m_i-1)/2\})=k=\mathrm{Card}(\{n_i=(m_i+1)/2\}).$$
Thus
$-\mathrm{deg}(f_*\mathcal{O}(m_1D_1+...+m_kD_k)/\mathcal{O}_C)/\mathrm{deg}\mathcal{L}$
is greater than or equal to
\begin{align*}
 &\underset{m_i\mathrm{\ even}}{\sum}\frac{m_i(m_i+2)}{8(m_i+1)}+\overset{}{\underset{n_i=(m_i-1)/2}{\sum}}\frac{(m_i-1)(m_i+1)}{8(m_i+1)}+\overset{}{\underset{n_i=(m_i+1)/2}{\sum}}\frac{(m_i+1)(m_i+3)}{8(m_i+1)}\\
&=\overset{}{\underset{m_i\mathrm{\ even}}{\sum}}\Big(\frac{m_i}{8}+\frac{m_i}{8(m_i+1)}\Big)+\overset{}{\underset{n_i=(m_i-1)/2}{\sum}}\frac{m_i-1}{8}+\overset{}{\underset{n_i=(m_i+1)/2}{\sum}}\frac{m_i+3}{8}\\
 &=\sum\frac{m_i}{8}+\underset{m_i\mathrm{\ even}}{\sum}\frac{m_i}{8(m_i+1)}+\frac{(-k)}{8}+\frac{3k}{8}\\
&=\frac{2g}{8}+\frac{1}{8}\Big(-2+\underset{m_i \mathrm{\
even}}{\sum}\frac{m_i}{m_i+1}+\underset{m_i\mathrm{\
odd}}{\sum}1\Big),
\end{align*}
so we get
$$L(C)=\frac{\mathrm{deg}
f_*\mathcal{O}(m_1D_1+...+m_kD_k)}{\mathrm{deg}(\mathcal{L})}+g$$
$$\qquad\qquad \qquad \leq \frac{3g}{4}-\frac{1}{8}\Big(-2+\underset{m_i\mathrm{\
even}}{\sum}\frac{m_i}{m_i+1}+\underset{m_i\mathrm{\
odd}}{\sum}1\Big).$$
\end{proof}
\begin{Remark}  D.Chen and M.M\"{o}ller \cite{Mo12} have obtained
a bound by using Cornalba-Harris-Xiao's slope inequality (the first
inequality) \cite[Theorem 2]{Xi}
$$L(C) \leq \frac{3g}{(g-1)}\kappa_{\mu}=\frac{g}{4(g-1)}\overset{k}{\underset{i=1}{\sum}}\frac{m_i(m_i+2)}{m_i+1}\leq \frac{3g}{4}.$$
\end{Remark}

\section{Weierstrass exponents}
This section is devoted  to the construction of the
Harder-Narasimhan filtration of $f_*\mathcal{O}(m_1D_1+...+m_kD_k)$
and the definition of Weierstrass exponents under some additional
assumptions.

\subsection{Weierstrass exponents}
If there is a filtration
\begin{equation}\label{exponent}
0\subset V_1\subset V_2 \subset... \subset
V_g=f_*\mathcal{O}(m_1D_1+...+m_kD_k)
\end{equation}
satisfying:
 (1)  $V_i/V_{i-1}$ is a line bundle, and
 (2)  $\mathrm{deg}(V_i/V_{i-1})$ is decreasing in $i$,
then it is the Harder-Narasimhan filtration of
$f_*\mathcal{O}(m_1D_1+...+m_kD_k)$, because each graded quotient
$V_i/V_{i-1}$ is already stable as it is a line bundle. The factors
$\mathrm{grad}_j(gr^{HN}_i)$ of  the Jordan-H\"{o}lder filtration of each
semistable graded quotient $\mathrm{gr}^{HN}_i$ are line bundles.

 We also get a filtration for
$f_*\omega_{S/C}$:
$$0\subset \mathcal{L}\otimes V_1\subset \mathcal{L}\otimes V_2\subset... \subset \mathcal{L}\otimes V_g=\mathcal{L} \otimes f_*\mathcal{O}(\sum m_iD_i)=f_*\omega_{S/C}.$$
\begin{Definition}
 If there exists a filtration as in \eqref{exponent}, we define the $i$-th \emph{Weierstrass exponent} $w_i$ as follows:
$$w_i=\frac{ \mathrm{deg}(\mathcal{L}\otimes V_i/\mathcal{L}\otimes V_{i-1})}{\mathrm{deg}(\mathcal{L})}.$$
\end{Definition}
\begin{Remark}It is obvious by definition that the sum of Weierstrass
exponents equals the sum of Lyapunov exponents
$$\overset{g}{\underset{i=1}{\sum}}\frac{ \mathrm{deg}(\mathcal{L}\otimes V_i/\mathcal{L}\otimes V_{i-1})}{\mathrm{deg}(\mathcal{L})}=\overset{g}{\underset{i=1}{\sum}} \frac{(\mathrm{deg}(\mathcal{L}\otimes V_i)-\mathrm{deg}(\mathcal{L}\otimes V_{i-1}))}{\mathrm{deg}(\mathcal{L})}=\frac{\mathrm{deg}(f_*\omega_{S/C})}{\mathrm{deg}(\mathcal{L})}.$$
\end{Remark}
When $H_{p_1,...,p_k}$ is nonvarying (i.e. $H_{p_1,...,p_k}\in
\cup_i WS_i$), there are many Weierstrass semigroup filtrations, and
we can construct the Harder-Narasimhan filtration recursively.
\begin{Theorem}
Assume that the Weierstrass semigroup $H_{p_1,...,p_k}$ is
nonvarying, then we can construct a filtration
$$0\subset V_1\subset V_2\subset... \subset V_g=f_*\mathcal{O}(m_1D_1+...+m_kD_k)$$
satisfying:

 1. $V_i/V_{i-1}$ is a line bundle, and

 2. $\mathrm{deg}(V_i/V_{i-1})$ is decreasing in $i$.
\end{Theorem}
\begin{proof}
For every  Weierstrass semigroup element $(d_1,...,d_k)$ of a
general fiber, we define the length of the element to be
$$l(d_1,...,d_k)=\mathrm{min}\{-d_1/(m_1+1),...,-d_k/(m_k+1)\},$$
where the fraction $-d_j/(m_j+1)$ is equal to
$\mathrm{deg}\mathcal{O}_{D_j}(d_jD_j)/\mathrm{deg}(\mathcal{L})$.

 For any vector bundle of the form
$f_*\mathcal{O}(a_1D_1+...+a_kD_k)$, we define the following set
$$L(a_1,...,a_k):=\{\underline{d}\in H_{p_1,...,p_k}|f_*\mathcal{O}(d_1D_1+...+d_kD_k)=f_*\mathcal{O}(a_1D_1+...+a_kD_k),\underline{d}\leq \underline{a}\}.$$
It is not empty because it contains the element $(d_1,...,d_k)\in
H_{p_1,...,p_k}$, for which the sum
$\overset{k}{\underset{i=1}{\sum}} d_i$ reaches the minumum when
$(d_1,...,d_k)$ varies in
$$\{\underline{d}|f_*\mathcal{O}(d_1D_1+...+d_kD_k)=f_*\mathcal{O}(a_1D_1+...+a_kD_k),\underline{d}\leq
\underline{a}\}.$$

We then construct the set $L_i$ and define the number $l_i$
recursively:
 $$L_g=
\{(m_1,...,m_k)\},\  l_g=l(m_1,...,m_k),$$
 $$......$$
$$L_i=
\{L(d_1,...,d_j-1,...,d_k)|(d_1,...,d_k)\in
L_{i+1},-d_j/(m_j+1)=l_{i+1}\},$$
$$l_i=\mathrm{min}\{l(d_1,...,d_k)|(d_1,...,d_k)\in L_i\},$$
 $$......$$

If $L_i\neq\emptyset$, then $l_i$ is defined, hence $L_{i-1}\neq\emptyset$ because $L(d_1,...,d_j-1,...,d_k)\neq\emptyset$  $(i\geq 2)$. So the definition makes sense.

It is obvious  that
$\mathrm{rk}(f_*\mathcal{O}(d_1D_1+...+d_kD_k))=i$ for any
$(d_1,...,d_k)\in L_i$.

For any $(e_1,...,e_k)\in L_{i-1}$, by our construction, there is a
$$(d_1,...,d_k)\in L_i,-d_j/(m_j+1)=l_i$$ 
such that $ (e_1,...,e_k)$
lies in $L(d_1,...,d_j-1,...,d_k)$. If we repeat the process from
$L_1$, then inductively we obtain a filtration
$$0\subset V_1\subset V_2\subset... \subset
V_g=f_*\mathcal{O}(m_1D_1+...+m_kD_k)$$ with $V_i/V_{i-1}$ being a
line bundle. The filtration is not unique because there maybe many
choices in each step.

 From the equalities
$V_i=f_*\mathcal{O}(d_1D_1+...+d_kD_k),$
$$V_{i-1}=f_*\mathcal{O}\Big(\sum e_iD_i\Big)=f_*\mathcal{O}(d_1D_1+...+(d_j-1)D_j+...+d_kD_k)$$
and the exact sequence
$$0\rightarrow V_{i-1}\rightarrow V_i\rightarrow\mathcal{O}_{D_j}(d_jD_j),$$
we get $$\mathrm{deg}(V_i/V_{i-1})/\mathrm{deg}(\mathcal{L})\leq
deg\mathcal{O}_{D_j}(d_jD_j)/\mathrm{deg}(\mathcal{L})=-d_j/(m_j+1)=l_i,$$
and
\begin{align*}
l_{i-1} &\geq min\{-e_1/(m_1+1),...,-e_k/(m_k+1)\}\\
&\geq
min\{-d_1/(m_1+1),...,-(d_j-1)/(m_j+1),...,-d_k/(m_k+1)\}\\
 &\geq
min\{-d_1/(m_1+1),...,-d_j/(m_j+1),...,-d_k/(m_k+1)\}\\
&=-d_j/(m_j+1)=l_i.
\end{align*}

When we assume that the Weierstrass semigroup is nonvarying, we get
$$\mathrm{deg}(V_i/V_{i-1})/\mathrm{deg}(\mathcal{L})=l_i$$ by Lemma
\ref{II}. Therefore, $\mathrm{deg}(V_i/V_{i-1})$ is decreasing in
$i$.
\end{proof}

\begin{Remark}From the proof, we can see that the Harder-Narasimhan filtration is constructed under
the weak assumption that   a subset of the Weierstrass semigroup is
nonvarying. This fact will be verified in hyperelliptic loci and
nonvarying strata of genus at most five.
\end{Remark}

\subsection{Hyperelliptic loci}
The square of any holomorphic 1-form $\omega$ on a hyperelliptic
curve fiber $F$ is a pullback $(\omega)^2=p^*q$ of some meromorphic
quadratic differential $q$ with simple poles   on $\mathbb{P}^1$
where the projection $p:F\rightarrow \mathbb{P}^1$ is the quotient
over the hyperelliptic involution. Denote by
$\mathcal{Q}(d_1,...,d_n)$ the stratification by orders of zeros and
simple poles of the corresponding quadratic differentials (see
\cite[section 2.2]{EKZ} for more details).

Let $F$ be the covering flat surface belonging to the stratum
 $\overline{\Omega\mathcal{M}}_g(m_1,...,m_k)$, then the
resulting holomorphic 1-form $\omega$ on  $F$ has zeros of the
following degrees:
\begin{itemize}
\item A zero $p_j$ of order $d$, of a meromorphic quadratic differential
$q$ on $\mathbb{P}^1$ gives rise to zeros on $F$ (\cite[p.
12]{EKZ}):

\item 1) Two zeros $p_{1j},p_{2j}$ of degree
 $m=d/2$,  when $d$ is even. In this case, $p_{1j},\ p_{2j}$ is a $g^1_2$,
 i.e. $h^0(p_{1j}+p_{2j})=2$.

\item 2) One zero $q_j$ of degree $m=d+1$, when $d>0$ is odd. In this
case, $q_j$ is a Weierstrass point, i.e. $h^0(2q_j)=2$.
\end{itemize}

Denote by $Q_j$ resp. $P_{1j}$ resp. $P_{2j}$ the section containing
$q_j$ resp. $p_{1j}$ resp. $p_{2j}$.

\begin{Proposition}  \label{hyperfiltration} Let $C$ be a Teichm\"{u}ller curve in the hyperelliptic locus of some
stratum $\overline{\Omega\mathcal{M}}_g(m_1,...,m_k)$, and denote by
$(d_1,...,d_n)$ the orders of singularities of the corresponding
quadratic differentials. Then there exists  a filtration
$$0\subset V_1\subset V_2\subset... \subset V_g=f_*\mathcal{O}(m_1D_1+...+m_kD_k)$$
satisfying:
 (1)  $V_i/V_{i-1}$ is a line bundle for each $i$,
 (2)  $\mathrm{deg}(V_i/V_{i-1})$ is decreasing in $i$.
\end{Proposition}
\begin{proof}

For each fiber $F$, the Weierstrass semigroup has at least a
nonvarying subgroup generated by the elements
$\{\{2q_j\}_{d_j\mathrm{\ odd}},\{p_{1j}+p_{2j}\}_{d_j \mathrm{\
even}}\}$, which is equal to:
$$\{\underset{d_j\mathrm{\ odd}}{\sum}2k_jq_j+\underset{d_j
\mathrm{\ even}}{\sum} n_j(p_{1j}+p_{2j})|2k_j\leq d_j+1,2n_j\leq
d_j+1\}$$ ($n_j\leq d_j/2  \Leftrightarrow  2n_j\leq d_j+1$ for
$d_j$ even).

 We order the following $g-1$ numbers
$\{\{-\frac{2k}{d_j+2}\}_{2k\leq d_j+1}\}$ to $\{N_1,...,N_{g-1}\}$
in decreasing order.

We transform $\{N_1,...,N_{g-1}\}$ to a new symbol set with $g-1$
element $\{T_1,...,T_{g-1}\}$ by using the following rule:
\begin{itemize}
\item when $d_j$ is odd
$$-\frac{2k}{d_j+2}\rightarrow 2Q_j,$$

\item when $d_j$ is even
$$ -\frac{2k}{d_j+2} \rightarrow
(P_{1j}+P_{2j}).$$
\end{itemize}

Then let $V_i=f_*\mathcal{O}(T_1+...+T_{i-1})$. We get a filtration
$$0\subset V_1\subset V_2\subset... \subset V_g=f_*\mathcal{O}(m_1D_1+...+m_kD_k).$$
If $N_i=-\frac{2k}{d_j+2}$ and $d_j$ is odd, then $\{N_1,...,N_i\}$
 contains $-\frac{2l}{d_j+2},l\leq k$, hence $T_1+...+T_{i-1}$ contain $2kQ_j$.

Since
$V_{i-1}=f_*\mathcal{O}(T_1+...+T_{i-1})=f_*\mathcal{O}(T_1+...+T_{i-1}+Q_j)$,
by the nonvarying property and Lemma \ref{II}, we get
$V_i/V_{i-1}=\mathcal{O}_{Q_j}(2kQ_j)$,
$$\mathrm{deg}(V_i/V_{i-1})=\mathrm{deg}(\mathcal{O}_{Q_j}(2kQ_j))=-\frac{2k}{(d_j+1)+1}=N_i.$$
If $N_i=-\frac{2k}{d_j+2}$ and $d_j$ is even, then $\{N_1,...,N_i\}$
just contains $-\frac{2l}{d_j+2},\ l\leq k$. $T_1+...+T_{i-1}$ just
contains $k(P_{1j}+P_{2j})$.

Similarly, we have
$V_{i-1}=f_*\mathcal{O}(T_1+...+T_{i-1})=f_*\mathcal{O}(T_1+...+T_{i-1}+P_{1j})$,
and  by the nonvarying property and Lemma \ref{II}, we get
$V_i/V_{i-1}=\mathcal{O}_{P_{2j}}(k(P_{1j}+P_{2j}))$, and that
$$\mathrm{deg}(V_i/V_{i-1})=\mathrm{deg}(\mathcal{O}_{P_{2j}}(k(P_{1j}+P_{2j})))=-\frac{k}{d_j/2+1}=N_i.$$
So $\mathrm{deg}(V_i/V_{i-1})$ is decreasing in $i$.
\end{proof}
\begin{Theorem} \label{hyperelliptic}Let $C$ be a Teichm\"{u}ller curve in the hyperelliptic locus of some
stratum $\overline{\Omega\mathcal{M}}_g(m_1,...,m_k)$, and denote by
$(d_1,...,d_n)$ the orders of singularities of the underlying
quadratic differentials. Then the Weierstrass exponents $w_i$ for
$C$ is the $i$-th largest number in the following set
$$\{1\}\cup\Big\{1-\frac{2k}{d_j+2}\Big\}_{\forall d_j, 0<2k\leq d_j+1}.$$
\end{Theorem}
\begin{proof}The Weierstrass exponents are as follows:
$$w_1=1,\ w_i=\mathrm{deg}(\mathcal{L}\otimes V_i/\mathcal{L}\otimes V_{i-1})/\mathrm{deg}(\mathcal{L})=1+N_{i-1}.$$
\end{proof}
\begin{Corollary}[{\cite[Corollary 1]{EKZ}}]\label{EKZ} Let $C$ be a Teichm\"{u}ller curve in the hyperelliptic locus of some stratum $\overline{\Omega\mathcal{M}}_g(m_1,...,m_k)$, and denote by
$(d_1,...,d_n)$ the orders of singularities of corresponding
quadratic differential. Then the sum of Lyapunov exponents of $C$ is
$$L(C)=\frac{1}{4}\sum_{d_j\mathrm{\ odd}} \frac{1}{d_j+2}.$$
\end{Corollary}
\begin{proof}Because $\overset{n}{\underset{i=1}{\sum}} d_i=-4$, we
have
 \begin{align*}
L(C)&=1+\underset{d_j \texttt{odd
}}{\sum}\overset{}{\underset{0<2k\leq
d_j+1}{\sum}}\Big(1-\frac{2k}{d_j+2}\Big)+\underset{d_j \texttt{even
}}{\sum}\overset{}{\underset{0<2k\leq
d_j+1}{\sum}}\Big(1-\frac{2k}{d_j+2}\Big)\\
&=1+\underset{0<d_j
\texttt{odd}}{\sum}\Big(\frac{d_j}{4}+\frac{1}{4(d_j+2)}\Big)+\underset{d_j
\texttt{even
}}{\sum}\frac{d_j}{4}=1+\underset{0<d_j \texttt{odd}}{\sum}\frac{1}{4(d_j+2)}+\underset{0<d_j}{\sum}\frac{d_j}{4}\\
 &=1+\underset{d_j
\texttt{odd}}{\sum}\frac{1}{4(d_j+2)}+\underset{d_j}{\sum}\frac{d_j}{4}=\frac{1}{4}\underset{d_j
\texttt{odd}}{\sum}\frac{1}{d_j+2}.
\end{align*}
\end{proof}

\subsection{The genus 3 case}
In what follows,  the dimension of special linear systems have been
discussed stratum by stratum in the corresponding section of
\cite[section 5, section 6, section 7]{CM11}. Note that for a
Teichm\"{u}ller curve, by Equation \eqref{VHS}, we have:
$$f_*\mathcal{O}(m_1D_1+...+m_kD_k)=\mathcal{O}_{C}\oplus
(f_*\mathcal{O}(m_1D_1+...+m_kD_k)/\mathcal{O}_{C}).$$ In the
stratum $\overline{\Omega\mathcal{M}}_3(3,1)$, a degenerate fibre is
not hyperelliptic, so $h^0(p_1+p_2)=h^0(2p_1)=1=g-2$ for any fiber.
By Lemma \ref{II}, we have
$$f_*\mathcal{O}(3D_1+D_2)=\mathcal{O}_{C}\oplus\mathcal{O}_{D_1}(3D_1)\oplus \mathcal{O}_{D_2}(D_2).$$

\begin{Remark}The Harder-Narasimhan filtration of
$f_*\mathcal{O}(3D_1+D_2)$ is
$$V_0=0\subset V_1=\mathcal{O}_{C} \subset V_2=\mathcal{O}_{C}\oplus\mathcal{O}_{D_2}(D_2)\subset
V_3=\mathcal{O}_{C}\oplus
\mathcal{O}_{D_2}(D_2)\oplus\mathcal{O}_{D_1}(3D_1).$$

The Harder-Narasimhan filtration of $f_*\omega_{S/C}$ is
$$0\subset \mathcal{L}\otimes V_1\subset \mathcal{L}\otimes V_2\subset \mathcal{L}\otimes V_3=\mathcal{L} \otimes f_*\mathcal{O}(3D_1+D_2).$$

The $i$-th Weierstrass exponent $w_i$ can be computed as follows:
\begin{align*}
    w_1 &=\mathrm{deg}(\mathcal{L}\otimes V_1/\mathcal{L}\otimes V_0)/\mathrm{deg}(\mathcal{L})\\
     &=\mathrm{deg}(\mathcal{L}\otimes \mathcal{O}_C)/\mathrm{deg}(\mathcal{L})=1,
\end{align*}
\begin{align*}
    w_2 & =\mathrm{deg}(\mathcal{L}\otimes V_2/\mathcal{L}\otimes V_1)/\mathrm{deg}(\mathcal{L})\\
     & =\mathrm{deg}(\mathcal{L}\otimes \mathcal{O}_{D_2}(D_2))/\mathrm{deg}(\mathcal{L})\\
     & =1+D_2^2/\mathrm{deg}(\mathcal{L})=\frac{1}{2},
\end{align*}
\begin{align*}
    w_3 &=\mathrm{deg}(\mathcal{L}\otimes V_3/\mathcal{L}\otimes V_2)/\mathrm{deg}(\mathcal{L}\\
     &=\mathrm{deg}(\mathcal{L}\otimes \mathcal{O}_{D_1}(3D_1))/\mathrm{deg}(\mathcal{L})\\
     &=1+3D_1^2/\mathrm{deg}(\mathcal{L})=\frac{1}{4},
\end{align*}
and the sum of Weierstrass exponents equals the sum of Lyapunov
exponents:
$$w_1+w_2+w_3=7/4=\lambda_1+\lambda_2+\lambda_3.$$
\end{Remark}

In the stratum $\overline{\Omega\mathcal{M}}^{odd}_3(2,2)$, the
theta characteristic is odd, so $h^0(p_1+p_2)=1=g-2$. By Lemma
\ref{II}, we have
$$f_*\mathcal{O}(2D_1+2D_2)=\mathcal{O}_{C}\oplus\mathcal{O}_{D_1}(2D_1)\oplus \mathcal{O}_{D_2}(2D_2).$$

In the stratum $\overline{\Omega\mathcal{M}}_3(2,1,1)$, we have
$h^0(p_1+p_3)=1=g-2$. By Lemma \ref{II}, we have
$$f_*\mathcal{O}(2D_1+D_2+D_3)=\mathcal{O}_{C}\oplus\mathcal{O}_{D_1}(2D_1)\oplus \mathcal{O}_{D_2}(D_2).$$
\begin{table}
\caption{genus 3}
\begin{tabular}{|c|c|c|c|c|c|}
  \hline
  zeros & component &  \multicolumn{3}{|c|}{Weierstrass exponents}  \\ \cline{3-5}
   &  & $w_2$& $w_3$ & $\sum w_i$  \\
   \hline
  (4)& hyp & 3/5 & 1/5 & 9/5  \\
    \hline
  (4) & odd & 2/5 & 1/5 & 8/5  \\
  \hline
   (3,1) &   & 2/4 & 1/4 & 7/4  \\
  \hline
   (2,2) & hyp & 2/3 & 1/3 & 2  \\
  \hline
   (2,2) & odd & 1/3 & 1/3 & 5/3  \\
  \hline
   (2,1,1) &   & 1/2 & 1/3 & 11/6  \\
  \hline
  (1,1,1,1) &   &   &   &   varying \\
  \hline
\end{tabular}
\end{table}
\subsection{The genus 4 case}
In the stratum $\overline{\Omega\mathcal{M}}_4(5,1)$, the curve $C$
does not meet the pointed Brill-Noether divisor $BN^1_{3,(2)}$, so
$h^0(3p_1)=1=g-3$.  By Lemma \ref{II},
$$f_*\mathcal{O}(5D_1+D_2)=\mathcal{O}_{C}\oplus f_*\mathcal{O}_{2D_1}(5D_1)\oplus \mathcal{O}_{D_2}(D_2).$$
We also have $h^0(4p_1)-1=h^0(3p_1)=h^0(2p_1)=1$, since by Lemma
\ref{I}
$$f_*\mathcal{O}_{2D_1}(5D_1)=\mathcal{O}_{D_1}(5D_1)\oplus \mathcal{O}_{D_1}(4D_1).$$

In the stratum $\overline{\Omega\mathcal{M}}^{even}_4(4,2)$, we have
$h^0(2p_1+p_2)=2,h^0(3p_1+2p_2)=3$. If  $h^0(3p_1+p_2)=3$, then by
Riemann-Roch $h^0(p_1+p_2)=2$, hence $p_1$ and $p_2$ are in the same
component of the fiber $F$. This component admits an involution that
acts on the zeros of  $\omega$. But $p_1$ and $p_2$ have different
orders of zeros, so
they cannot be conjugate under the involution. This contradiction implies that we have $h^0(3p_1+p_2)=2$ and $h^0(p_1+p_2)=1$.\\
We get a Weierstrass pair filtration
$$\{\{(1,1),(2,1)\},\{(3,1),(3,2)\},\{(3,2),(4,2)\}\},$$
that is
$$0\subset f_*\mathcal{O}(D_1+D_2)\subset
f_*\mathcal{O}(2D_1+D_2)=f_*\mathcal{O}(3D_1+D_2)$$  $$\subset
f_*\mathcal{O}(3D_1+2D_2)=f_*\mathcal{O}(3D_1+2D_2)\subset
f_*\mathcal{O}(4D_1+2D_2).$$

 This is also
the Harder-Narasimhan filtration, as the graded quotients are line
bundles.

In the stratum $\overline{\Omega\mathcal{M}}^{odd}_4(4,2)$,  the
theta characteristic is odd, so $h^0(2p_1+p_2)=1=g-3$. By Lemma
\ref{II},
 $$f_*\mathcal{O}(4D_1+2D_2)=\mathcal{O}_{C}\oplus f_*\mathcal{O}_{2D_1}(4D_1)\oplus  \mathcal{O}_{D_2}(2D_2).$$
We also have $h^0(3p_1+p_2)-1=h^0(2p_1+p_2)=h^0(p_1+p_2)=1$, since
by Lemma \ref{I}
$$f_*\mathcal{O}_{2D_1}(4D_1)=\mathcal{O}_{D_1}(4D_1)\oplus \mathcal{O}_{D_1}(3D_1).$$

In the stratum $\overline{\Omega\mathcal{M}}^{non-hyp}_4(3,3)$, the
curve $C$ does not meet the pointed Brill-Noether divisor
$BN^1_{3,(1,1)}$, so $h^0(2p_1+p_2)=1=g-3$. By Lemma \ref{II},
$$f_*\mathcal{O}(3D_1+3D_2)=\mathcal{O}_{C}\oplus f_*\mathcal{O}_{D_1}(3D_1)\oplus  \mathcal{O}_{2D_2}(3D_2).$$
We have $h^0(3p_1+p_2)-1=h^0(2p_1+p_2)=h^0(p_1+p_2)=1$, since by
Lemma \ref{I}
$$f_*\mathcal{O}_{2D_2}(3D_2)=\mathcal{O}_{D_2}(3D_2)\oplus \mathcal{O}_{D_2}(2D_2).$$

In the stratum $\overline{\Omega\mathcal{M}}^{odd}_4(2,2,2)$,
 by Clifford's Theorem, we get $h^0(p_1+p_2+p_3)=1$. By Lemma \ref{II} we get
$$f_*\mathcal{O}(2D_1+2D_2+2D_3)=\mathcal{O}_{C}\oplus\mathcal{O}_{D_1}(2D_1)\oplus \mathcal{O}_{D_2}(2D_2)\oplus \mathcal{O}_{D_3}(2D_3).$$

In the stratum  $\overline{\Omega\mathcal{M}}_4(3,2,1)$,
 the curve $C$ does not
meet the pointed Brill-Noether divisor $BN^1_{4,(1,1,2)}$, as it has
been shown in \cite[section 6.8]{CM11}  that
$h^0(2p_1+p_2)=h^0(p_1+p_2+p_3)=1=g-3$. By Lemma \ref{II},
$$f_*\mathcal{O}(3D_1+2D_2+D_3)=\mathcal{O}_{C}\oplus f_*\mathcal{O}_{2D_1}(3D_1)\oplus \mathcal{O}_{D_2}(2D_2).$$
We have $h^0(3p_1+p_2)-1=h^0(2p_1+p_2)=h^0(p_1+p_2)=1$, since by
Lemma \ref{I}
$$f_*\mathcal{O}_{2D_1}(3D_1)=\mathcal{O}_{D_1}(3D_1)\oplus \mathcal{O}_{D_1}(2D_1).$$

\begin{table}
\caption{genus 4}
 \begin{tabular}{|c|c|c|c|c|c|c|}
  \hline
  zeros & component &  \multicolumn{4}{|c|}{Weierstrass exponents}\\ \cline{3-6}
   &  & $w_2$& $w_3$ &$w_4$ & $\sum w_i$  \\
   \hline
  (6)&  hyp & 5/7 & 3/7 & 1/7& 16/7 \\
    \hline
  (6) & even& 4/7 & 2/7 & 1/7 &14/7   \\
  \hline
  (6) & odd & 3/7 & 2/7 & 1/7 &13/7   \\
  \hline
   (5,1) &  & 1/2 & 2/6 & 1/6 & 2   \\
   \hline
   (3,3) & hyp & 3/4 & 2/4 & 1/4 & 5/2  \\
   \hline
   (3,3) & non-hyp & 2/4 & 1/4 & 1/4 & 2  \\
  \hline
     (4,2) & even & 3/5 & 1/3 & 1/5 &32/15   \\
  \hline
   (4,2) & odd & 2/5 & 1/3 & 1/5 &29/15   \\
  \hline
  (2,2,2) & odd & 1/3 & 1/3 & 1/3 &2 \\
  \hline
     (3,2,1) &   & 1/2 & 1/3 & 1/4 &25/12   \\
  \hline
\end{tabular}
\end{table}
\subsection{The genus 5 case}
In the stratum  $\overline{\Omega\mathcal{M}}_5(5,3)$, the curve $C$
does not meet the pointed Brill-Noether divisor $BN^1_{4,(1,2)}$,
therefore $h^0(2p_1+2p_2)=1$, and by Riemann-Roch Theorem, we get
$h^0(3p_1+p_2)=1=g-4$. By Lemma \ref{II},
$$f_*\mathcal{O}(5D_1+3D_2)=\mathcal{O}_{C}\oplus\mathcal{O}_{2D_1}(5D_1)\oplus \mathcal{O}_{2D_2}(3D_2).$$
Thus $h^0(4p_1+p_2)-1=h^0(3p_1+p_2)=h^0(2p_1+p_2)=1$, and
$h^0(p_1+4p_2)-1=h^0(p_1+3p_2)=h^0(p_1+2p_2)=1$, so by Lemma \ref{I}
$$f_*\mathcal{O}_{2D_1}(5D_1)=\mathcal{O}_{D_1}(5D_1)\oplus \mathcal{O}_{D_1}(4D_1),$$
$$f_*\mathcal{O}_{2D_2}(3D_2)=\mathcal{O}_{D_2}(3D_2)\oplus \mathcal{O}_{D_2}(2D_2).$$

In the stratum  $\overline{\Omega\mathcal{M}}^{odd}_5(6,2)$, we need
the following lemma, whose proof is due to D.Chen. The reader can
also read \cite[Lemma A.8]{CM12} for another proof.
\begin{Lemma}[Chen]\label{chen}In the stratum
$\overline{\Omega\mathcal{M}}^{odd}_5(6,2)$, we have
$h^0(3p_1+p_2)=1$.
\end{Lemma}
 \begin{proof}If the curve $X$ is smooth or irreducible, then
$h^0(3p_1+p_2)=3$ would imply, by Clifford's Theorem, that $X$ is
hyperelliptic and $2p_1$ linear equivalence to $ p_1 + p_2$ is
impossible.

If $X$ is reducible, it can have at most two components $Z$ and $Y$
meeting at $n$ nodes ($n>1$), such that $6=2g_1-2+n$ and
$2=2g_2-2+n$. Therefore the only possibilities for $(g_1,g_2,n)$ are
$(3,1,2)$ and $(2,0,4)$.

For the first case, the elliptic component $Y$ contains $p_2$ and
$h^0(p_2)=1$, i.e. all the sections are given by constant functions.
The other component $Z$ contains $p_1$ with $h^0(3p_1)<3$, and a
section on $Z$ uniquely determines the constant section on Y, by its
values at the nodes (assuming the same value). Hence
$h^0(3p_1+p_2)<3$ on $X$.

For the second case, $h^0(p_2)=2$ on the rational component $Y$.
Then $h^0(3p_1)$ has to be $2$ on $Z$, hence $p_1$ is a Weierstrass
point. But in order to glue two sections on $Y$ and $Z$, they need
to have the same value at each of the four nodes. The four nodes
form two conjugate pairs on the hyperelliptic curve $Z$, hence
gluing two sections still imposes two conditions. Therefore
$h^0(3p_1+p_2)\leq 2+2-2<3$ on $X$.
\end{proof}

We get $h^0(3p_1+p_2)=1=g-4$. By Lemma \ref{II},
$$f_*\mathcal{O}(6D_1+2D_2)=\mathcal{O}_{C}\oplus\mathcal{O}_{3D_1}(6D_1)\oplus \mathcal{O}_{D_2}(2D_2),$$
and we have $$h^0(4p_1+p_2)-1=h^0(3p_1+p_2)=h^0(2p_1+p_2)=1,$$ since
by Lemma \ref{I}
$$f_*\mathcal{O}_{3D_1}(6D_1)=\mathcal{O}_{D_1}(6D_1)\oplus \mathcal{O}_{D_1}(5D_1)\oplus \mathcal{O}_{D_1}(4D_1).$$

\begin{table}
\caption{genus 5}
\begin{tabular}{|c|c|c|c|c|c|c|c|}
  \hline
  zeros & component &  \multicolumn{5}{|c|}{Weierstrass exponents} \\ \cline{3-7}
   &  & $w_2$& $w_3$ &$w_4$ &$w_5$ & $\sum w_i$  \\
   \hline
  (8)&  hyp & 7/9 & 5/9 & 3/9 &1/9& 25/9 \\
    \hline
  (8) & even& 5/9 & 3/9 & 2/9 &1/9&20/9   \\
  \hline
  (8) & odd & 4/9 & 3/9 & 2/9 &1/9&19/9   \\
  \hline
   (5,3) &   & 1/2 & 1/3 & 1/4 &1/6& 9/4   \\
  \hline
      (6,2) & odd  & 3/7 & 1/3 & 2/7 &1/7& 46/21  \\
  \hline
     (4,4) & hyp  & 4/5 & 3/5 & 2/5 &1/5& 3  \\
  \hline
\end{tabular}
\end{table}
\begin{Theorem}\label{lowgenus}For a Teichm\"{u}ller curve in the strata\\
$\overline{\Omega\mathcal{M}}_3(3,1),\overline{\Omega\mathcal{M}}^{odd}_3(2,2),\overline{\Omega\mathcal{M}}_3(2,1,1)$\\
$\overline{\Omega\mathcal{M}}_4(5,1),\overline{\Omega\mathcal{M}}^{odd}_4(4,2),\overline{\Omega\mathcal{M}}^{non-hyp}_4(3,3),\overline{\Omega\mathcal{M}}^{odd}_4(2,2,2),\overline{\Omega\mathcal{M}}_4(3,2,1)$\\
$\overline{\Omega\mathcal{M}}_5(5,3),\overline{\Omega\mathcal{M}}^{odd}_5(6,2),$\\
the Weierstrass exponents can explicitly be calculated as in Tables
1, 2 and 3. Moreover $f_*{\omega_{S/C}}$ splits into a direct sum of
line bundles.

For a Teichm\"{u}ller curve in the stratum
$\overline{\Omega\mathcal{M}}^{even}_4(4,2)$, the Weierstrass
exponents can explicitly be calculated as in Table 2.
\end{Theorem}
\begin{proof} In all the cases, we have constructed a filtration
of type \eqref{exponent}.
\end{proof}
\section*{Acknowledgement}We are very grateful to Dawei Chen
for pointing out a mistake in a preliminary version of this paper
and for making several remarks. We thank Martin M\"{o}ller and Anton
Zorich for their helpful advice. We also thank Ke Chen and Ronggang
Shi for carefully reading the paper and their suggestions on the
presentation of the paper. We would like to thank anonymous referees
for their careful reviews of the paper that include many
important points.

\end{document}